\documentclass[twoside,11pt]{article}

\usepackage{amssymb, amsmath}

\newcommand{\rk}{{\rm rk}}

\newcommand{\st}{{\rm st}}
\newcommand{\tot}{{\rm tot}}

\newcommand{\depth}{{\rm depth} \,}

\newcommand{\injdim}{{\rm injdim}}
\newcommand{\GKdim}{{\rm GKdim}}
\newcommand{\Cdim}{{\rm Cdim}}
\newcommand{\Hom}{{\rm Hom}}
\newcommand{\HOM}{{\rm \underline{\Hom}}}
\newcommand{\Ext}{{\rm Ext}}
\newcommand{\EXT}{{\rm \underline{\Ext}}}

\newcommand{\GrMod}{{\sf GrMod}}

\newcommand{\tq}{\,|\,}
\newcommand{\proof}{\noindent {\it Proof.$\;$}}
\newcommand{\qed}{\hfill \rule{1.5mm}{1.5mm}}

\newcommand\zum[2]{\sum_{\substack{#1\\#2}}}
\newcommand\zzum[3]{\sum_{\substack{#1\\#2\\#3}}}

\newtheorem{theorem}{Theorem}[section]

\newtheorem{remark}[theorem]{Remark}

\newtheorem{subtheorem}{Theorem}[subsection]
\newtheorem{subproposition}[subtheorem]{Proposition}
\newtheorem{subdefinition}[subtheorem]{Definition}
\newtheorem{sublemma}[subtheorem]{Lemma}

\newtheorem{subcorollary}[subtheorem]{Corollary}
\newtheorem{subremark}[subtheorem]{Remark}

\newcommand{\titre}{Quantum graded algebras with a straightening law
and the AS-Cohen-Macaulay property for quantum determinantal rings
and quantum grassmannians.}
\newcommand{\fichier}{T.H. Lenagan and L. Rigal}
\newcommand{\tit}{Quantum graded algebras with a straightening law.}

\begin{document}

\title{{\vspace{-1.5cm} \bf \titre}}
\date{T.H. Lenagan and L. Rigal 
\thanks{This research was supported by the Leverhulme Research 
Interchange Grant
F/00158/X.  Part of the work was done while the authors were visiting
Antwerp supported by the European Science Foundation programme
`Noncommutative Geometry', and while the second author was visiting the
University of Edinburgh supported by a grant from the Edinburgh
Mathematical Society.}}
\maketitle

\begin{abstract} 
 
We study quantum analogues of quotient varieties, namely quantum
grassmannians and quantum determinantal rings, from the point of view of
regularity conditions.  More precisely, we show that these rings are
AS-Cohen-Macaulay and determine which of them are AS-Gorenstein.  Our
method is inspired by the one developed by De Concini, Eisenbud and
Procesi in the commutative case.  Thus, we introduce and study the
notion of a {\em quantum graded algebra with a staightening law} on a
partially ordered set, showing in particular that, among such algebras,
those whose underlying poset is wonderful are AS-Cohen-Macaulay.  Then,
we prove that both quantum grassmannians and quantum determinantal rings
are quantum graded algebras with a staightening law on a wonderful
poset, hence showing that they are AS-Cohen-Macaulay.  In this last
step, we are lead to introduce and study (to some extent) natural
quantum analogues of Schubert varieties. 

\end{abstract}

\vskip .5cm
\noindent
{\em 2000 Mathematics subject classification:} 16W35, 16P40, 16S38, 17B37,
20G42.

\vskip .5cm
\noindent
{\em Key words:} Quantum matrices, quantum grassmannian, quantum Schubert 
variety, quantum determinantal ring, straightening laws, Cohen-Macaulay,
Gorenstein. 



\section*{Introduction.}
 
Let $\Bbbk$ be an arbitrary field, $m,n$ be positive integers such that 
$m \le n$ and $q$ be a non-zero element of $\Bbbk$. We denote by
${\mathcal O}_q(M_{m,n}(\Bbbk))$ the quantum deformation
of the coordinate ring of the variety $M_{m,n}(\Bbbk)$ of $m \times n$ matrices
with entries in $\Bbbk$ and by ${\mathcal O}_q(G_{m,n}(\Bbbk))$
the quantum deformation of the homogeneous coordinate ring of the 
Grassmannian $G_{m,n}(\Bbbk)$ of $m$-dimensional subspaces
of $\Bbbk^n$ (see Definitions \ref{def-q-mat} and \ref{def-q-grass}).
In addition, if $1 \le t \le m$, we denote by ${\mathcal I}_t$ the ideal of
${\mathcal O}_q(M_{m,n}(\Bbbk))$ generated by the $t \times t$ quantum minors
(see the beginning of Subsection \ref{subsec-q-grass}). 
As is well known, the algebras ${\mathcal O}_q(M_{m,n}(\Bbbk))$ and 
${\mathcal O}_q(G_{m,n}(\Bbbk))$ are noncommutative analogues of 
${\mathcal O}(M_{m,n}(\Bbbk))$
(the coordinate ring of the affine variety $M_{m,n}(\Bbbk)$) 
and of ${\mathcal O}(G_{m,n}(\Bbbk))$ (the homogeneous coordinate ring of
the projective variety $G_{m,n}(\Bbbk)$) 
in the sense that the usual coordinate rings of these 
varieties are recovered when the parameter $q$ is taken to be $1$
(and the base field is algebraically closed).
In the same way, ${\mathcal O}_q(M_{m,n}(\Bbbk))/{\mathcal I}_t$ is a quantum
deformation of the coordinate ring of the determinantal variety 
$M_{m,n}^{\le t-1}(\Bbbk)$ of those matrices in $M_{m,n}(\Bbbk)$ whose rank
is at most $t-1$.

The commutative rings ${\mathcal O}(M_{m,n}(\Bbbk))/{\mathcal I}_t$ and 
${\mathcal O}(G_{m,n}(\Bbbk))$ have been extensively studied in the past 
decades; an extensive account of their interest and properties can be
found, for example, in [BV]. They are rings of invariants for natural
actions of the general and special linear groups on a suitable 
polynomial ring. It follows that they are normal domains. In addition,
if $\Bbbk$ is an algebraically closed field of characteristic zero, 
it follows by the Hochster-Roberts theorem (see [BH; Theorem 6.5.1]) 
that these rings are Cohen-Macaulay. 
In fact, the Cohen-Macaulay property of these rings still
holds over an arbitrary field. One way to prove this is to make use of
the notion of {\em graded algebras with a staightening law}, initiated in works
by De Concini, Eisenbud and Procesi (see [E] and [DEP]). For a detailed 
description of this method and extensive comments on the historical 
background, the reader is referred to [BV] and [BH]. 

Hence, ${\mathcal O}_q(M_{m,n}(\Bbbk))/{\mathcal I}_t$ and
${\mathcal O}_q(G_{m,n}(\Bbbk))$ can be thought of as quantum analogues
of quotient varieties under the action of reductive groups. 
In support of 
this point of view, note that K.R. Goodearl, A.C. Kelly and the authors have 
proved that 
${\mathcal O}_q(M_{m,n}(\Bbbk))/{\mathcal I}_t$ 
is the ring of co-invariants for a natural
co-action of the Hopf algebra ${\mathcal O}_q(GL_t(\Bbbk))$ on 
${\mathcal O}_q(M_{m,t}(\Bbbk)) \otimes {\mathcal O}_q(M_{t,n}(\Bbbk))$,
see [GLR], and 
${\mathcal O}_q(G_{m,n}(\Bbbk))$ is the ring of co-invariants for a
natural
co-action of the Hopf algebra ${\mathcal O}_q(SL_m(\Bbbk))$ on
${\mathcal O}_q(M_{m,n}(\Bbbk))$, see  [KLR].
In this context, it is natural to study 
${\mathcal O}_q(M_{m,n}(\Bbbk))/{\mathcal I}_t$ and
${\mathcal O}_q(G_{m,n}(\Bbbk))$ from the point of view of  regularity
properties, where by regularity properties 
we mean the maximal order property,
the AS-Cohen-Macaulay property and the AS-Gorenstein property, which are 
noncommutative analogues of classical regularity properties in commutative 
algebra and algebraic geometry.
In [KLR] and [LR], the maximal order property for the above rings 
has been investigated. Recall that the maximal order property is a 
non commutative
analogue for normality.
Here, our aim is to study the $\Bbbk$-algebras
${\mathcal O}_q(G_{m,n}(\Bbbk))$ and
${\mathcal O}_q(M_{m,n}(\Bbbk))/{\mathcal I}_t$ from the point
of view of the AS-Cohen-Macaulay and AS-Gorenstein properties, as 
defined in recent works on
noncommutative algebraic geometry (see, for example, [JZ]). These properties are
noncommutative analogues of the notions of Cohen-Macaulay and Gorenstein
rings.

There are many problems associated with obtaining quantum
analogues of these famous commutative results. Of course, the most
obvious problems arise due to the lack of commutativity; this forces us to study
in detail specific kinds of commutation relations. Note
that, although the algebras we consider are known to be algebras of co-invariants
of suitable Hopf algebra coactions, at the moment we have no
known methods for exploiting this fact. \\

In order to prove that ${\mathcal O}_q(M_{m,n}(\Bbbk))/{\mathcal I}_t$
and ${\mathcal O}_q(G_{m,n}(\Bbbk))$ are AS-Cohen-Macaulay $\Bbbk$-algebras,
we introduce and study the notion of {\em quantum graded algebras with a 
staightening  law} (quantum graded A.S.L. for short). Roughly speaking,
a quantum graded A.S.L. is an ${\mathbb N}$-graded algebra with a 
partially ordered finite set $\Pi$ of homogeneous generators, 
satisfying the following properties: standard monomials (that is 
products of elements of $\Pi$ in increasing order) 
form a free family; the product of any two 
incomparable elements of $\Pi$ can be written as a linear combination
of standard monomials in a way compatible with the partial order on $\Pi$
(these are the so-called {\em straightening relations});
given two elements $\alpha,\beta\in\Pi$, there exists a skew commutator
of $\alpha$ and $\beta$ which can be written as a linear combination
of standard monomials in a way compatible with the partial order on $\Pi$ 
(these are the so-called {\em commutation relations}). 
We then show that quantum graded A.S.L. are AS-Cohen-Macaulay (under a mild 
assumption on $\Pi$). The main point then is to show that
${\mathcal O}_q(M_{m,n}(\Bbbk))/{\mathcal I}_t$ and 
${\mathcal O}_q(G_{m,n}(\Bbbk))$ are examples of such quantum graded
A.S.L.
We expect that the notion of a quantum graded A.S.L. will prove useful
in many other contexts. 
\\

The paper is organised as follows. In Section 1, we introduce and study
the notion of quantum graded A.S.L. on a partially ordered finite 
set $\Pi$. In particular, we show that such algebras are noetherian and we 
compute their Gelfand-Kirillov dimension. We also study in detail certain 
ideals of a quantum graded A.S.L., namely those which are generated by 
$\Pi$-ideals (see the beginning of Subsection \ref{subsect-pi-ideals}). 
In particular,  we show that these ideals are generated by normalising 
sequences.
The $\Pi$-ideals are essential for our arguments 
because they allow us to study quantum graded A.S.L. by
inductive arguments. In Section 2, we first introduce 
the homological background necessary to define and study the AS-Cohen-Macaulay 
and AS-Gorenstein properties and then
recall the results concerning these notions 
that we will need. 
Then we prove that a quantum graded A.S.L. is AS-Cohen-Macaulay, 
provided that the poset $\Pi$ is {\em wonderful}. 
The aim of Section 3 is then
to show that the algebras ${\mathcal O}_q(M_{m,n}(\Bbbk))/{\mathcal I}_t$
and ${\mathcal O}_q(G_{m,n}(\Bbbk))$ are quantum graded A.S.L. 
For this, we first study ${\mathcal O}_q(G_{m,n}(\Bbbk))$ in detail: this 
essentially consists of proving the existence of straightening relations
and specific kinds of commutation relations between maximal minors of 
${\mathcal O}_q(M_{m,n}(\Bbbk))$. These two types of relations are established 
by introducing and studying natural quantum analogues of (coordinate rings on)
Schubert varieties over arbitrary integral domains. In this way, we prove that
${\mathcal O}_q(G_{m,n}(\Bbbk))$ is a quantum graded A.S.L. and deduce that
${\mathcal O}_q(M_{m,n}(\Bbbk))$ enjoys the same property by means of a
dehomogenisation map that relates these two algebras. 
It is then easy to show that the $\Bbbk$-algebras
${\mathcal O}_q(M_{m,n}(\Bbbk))/{\mathcal I}_t$ are also quantum graded A.S.L.
Finally, in Section 4, we deduce that the $\Bbbk$-algebras
${\mathcal O}_q(M_{m,n}(\Bbbk))/{\mathcal I}_t$ and 
${\mathcal O}_q(G_{m,n}(\Bbbk))$ are AS-Cohen-Macaulay and determine which
of them are AS-Gorenstein by means of a criterion involving 
their Hilbert series.\\
  
If $S$ is a finite set, we denote its cardinality by $|S|$.\\

We would like to thank the many colleagues with whom we have discussed this work,
especially J. Alev, G. Cauchon, K.R. Goodearl, P. J\o rgensen, S. Launois, 
S.P. Smith and J.J. Zhang.

\section{Quantum graded algebras with a straightening law.}

Throughout this section, $\Bbbk$ is a field.  An ${\mathbb N}$-graded
$\Bbbk$-algebra is a $\Bbbk$-algebra together with a family of
$\Bbbk$-subspaces $A_i$, for $i \in {\mathbb N}$, such that $A=\oplus_{i
\in {\mathbb N}} A_i$ and $A_i A_j \subseteq A_{i+j}$, for each pair
$i,j \in {\mathbb N}$.  An ${\mathbb N}$-graded $\Bbbk$-algebra is said
to be {\it connected} if $A_0=\Bbbk .  1$.  An ${\mathbb N}$-graded
$\Bbbk$-algebra is connected if and only if it is generated as a
$\Bbbk$-algebra by $\oplus_{i \ge 1} A_i$. 

\subsection{Definition and elementary properties.}

Let $A$ be an ${\mathbb N}$-graded algebra and $\Pi$ a finite subset of 
elements of $A$ with a partial order $<_\st$. A  
{\em standard monomial} on $\Pi$ is an element
of $A$ which is either $1$ or of the form $\alpha_1\dots\alpha_s$, 
for some $s\geq 1$, where $\alpha_1,\dots,\alpha_s \in \Pi$ and
$\alpha_1\le_\st\dots\le_\st\alpha_s$. 

\begin{subdefinition} \label{q-gr-asl} -- 
Let $A$ be an ${\mathbb N}$-graded $\Bbbk$-algebra
and $\Pi$ a finite subset of $A$ equipped with a partial order $<_\st$. 
We say that $A$ is a {\em quantum graded algebra with a 
straightening law} 
({\em quantum graded A.S.L.} for short) on the poset $(\Pi,<_\st)$ 
if the following conditions are satisfied.\\ 

\noindent 
(1) The elements of $\Pi$ are homogeneous with positive degree.\\
(2) The elements of $\Pi$ generate $A$ as a $\Bbbk$-algebra.\\
(3) The set of standard monomials on $\Pi$ is a free family.\\
(4) If $\alpha,\beta\in\Pi$ are not comparable for $<_\st$, then $\alpha\beta$ 
is a linear combination of terms $\lambda$ or $\lambda\mu$, where 
$\lambda,\mu\in\Pi$, $\lambda\le_\st\mu$ and $\lambda<_\st\alpha,\beta$.\\
(5) For all $\alpha,\beta\in\Pi$, there exists $c_{\alpha\beta} \in \Bbbk^\ast$ 
such that $\alpha\beta-c_{\alpha\beta}\beta\alpha$ is a linear combination of 
terms $\lambda$ or $\lambda\mu$, where $\lambda,\mu\in\Pi$,
$\lambda\le_\st\mu$ and $\lambda<_\st\alpha,\beta$.
\end{subdefinition}

Notice that, in the above definition, the case $\Pi=\emptyset$ is not excluded.
Hence, $\Bbbk$ is a quantum graded A.S.L. on $\emptyset$ and, of course, it is
the only one. 

\begin{subremark} \label{rem-before-rank} -- 
\rm Let $A$ be a quantum graded A.S.L. on the set $\Pi$.
If $\omega$ is a standard monomial on $\Pi$ and $\omega\neq 1$, then 
condition (3) of
Definition \ref{q-gr-asl} implies that it can be written in a unique way as a 
product $\alpha_1\dots\alpha_s$, for some $s \geq 1$, where 
$\alpha_1,\dots,\alpha_s \in \Pi$ and $\alpha_1\le_\st\dots\le_\st\alpha_s$. 
In this case, we say that $\omega$
is a standard monomial of {\em length $s$}. For convenience, we say that $1 \in
A$ is the (unique) standard monomial of length $0$.
\end{subremark}

Let $\pi$ be an element of a finite poset $\Pi$ with partial order $<_\st$; 
we define the {\em rank} 
of $\pi$, denoted $\rk \pi$, as in [BV; Chap. 5.C, p. 55]. 
Thus, $\rk \pi=k$ if and only if there is a chain $\pi_1 <_\st \dots <_\st
\pi_{k-1} <_\st \pi_k=\pi$, with 
$\pi_i \in \Pi$ for $1 \le i \le k$ and no such 
chain of
greater length exists. If $\Omega$ is a subset of $\Pi$, we define its rank by
$\rk\Omega=\max\{\rk \pi,\,\pi\in\Omega\}$. The rank of the empty set is taken
to be $0$. \\

Let $A$ be a quantum graded A.S.L. on the set $\Pi$ partially ordered by 
$<_\st$. For $s \in {\mathbb N}^\ast$, we denote by $\Pi^s$ the set of elements 
of $A$ which either equal $1$ or can be written as a product of $t$ elements of 
$\Pi$, with 
$1 \le t \le s$. In addition, we let $V^s$ be the $\Bbbk$-subspace of 
$A$ generated by $\Pi^s$.
On the other hand, we denote by $\Pi^s_\st$ the set of those elements of $\Pi^s$
which are standard monomials of length at most equal to $s$ and by $V^s_\st$ the
$\Bbbk$-subspace of $A$ generated by $\Pi^s_\st$. For convenience, we put
$\Pi^0=\Pi^0_\st=\{1\}$ and $V^0=V^0_\st=\Bbbk$. Clearly, 
$V^s_\st\subseteq V^s$, for $s \in {\mathbb N}$.
Of course,
since $\Pi^1_\st=\Pi^1=\Pi\cup\{1\}$, we have $V^1_\st = V^1$. Also, as an easy
consequence of conditions (4) and (5) in the definition of a quantum graded 
A.S.L., we have $V^2_\st = V^2$.

\begin{sublemma} \label{pre-pre-st-fam-gen} -- 
Let $s \in {\mathbb N}$, $s\ge 2$ and assume that $V_\st^s=V^s$. 
If $\pi_1,\dots,\pi_{s+1}$ are elements of $\Pi$ such that 
$\pi_2 \le_\st\dots\le_\st\pi_{s+1}$, then $\pi=\pi_1\dots\pi_{s+1} 
\in V^{s+1}_\st$.
\end{sublemma}

\proof We proceed by induction on $\rk\pi_1$. If $\rk \pi_1=1$, then $\pi_1$ is
minimal in $\Pi$ with respect to $<_\st$. Two cases may occur. If $\pi_1$ and 
$\pi_2$ are not comparable with respect to $<_\st$, then the minimality of 
$\pi_1$ and condition (4) in the definition of a quantum graded A.S.L. imply 
that $\pi_1\pi_2=0$ and thus $\pi=0$. Otherwise, $\pi_1,\pi_2$ are comparable 
and, since $\pi_1$ is minimal, we have $\pi_1 \le_\st \pi_2$. It follows that 
$\pi$ is a standard monomial. In both cases, $\pi \in V_\st^{s+1}$. 

Now, assume that the result is true when $1 \le \rk\pi_1 \le r$, where $r
< \rk \Pi$ and consider an element $\pi$ as above, for which
$\rk\pi_1=r+1$.  Thus, $\pi=\pi_1\dots\pi_{s+1}$, where
$\pi_1,\pi_2,\dots,\pi_{s+1}$ are elements of $\Pi$ such that
$\pi_2\le_\st\dots\le_\st\pi_{s+1}$. 

First case: $\pi_1$ and $\pi_2$
are not comparable.  Condition (4) in the definition of a quantum graded
A.S.L.  shows that $\pi$ is a linear combination of terms of
$V^s=V_\st^s$ and of terms
$\pi_1^{\prime\prime}\dots\pi_{s+1}^{\prime\prime}$ such that
$\pi_1^{\prime\prime} <_\st \pi_1$.  

Second case: $\pi_1$ and $\pi_2$
are comparable.  If $\pi_1 \le_\st \pi_2$, $\pi$ is standard and thus in
$V_\st^{s+1}$.  If $\pi_1 >_\st \pi_2$, condition (5) in the definition
of a quantum graded A.S.L.  shows that $\pi$ is the sum of $\omega\in
V_\st^s$ and of a linear combination of $\pi_2\pi_1\pi_3\dots\pi_{s+1}$
together with terms of the form $\pi_1^\prime\dots\pi_{s+1}^\prime$ such
that $\pi_1^\prime <_\st \pi_1$.  Hence, all the terms in this linear
combination are elements of the form
$\pi_1^{\prime\prime}\dots\pi_{s+1}^{\prime\prime}$ with
$\pi_1^{\prime\prime} <_\st \pi_1$.  

Thus, in all cases,
$\pi_1\dots\pi_{s+1}$ is the sum of an element of $V_\st^{s+1}$ and of a
linear combination of terms
$\pi_1^{\prime\prime}\dots\pi_{s+1}^{\prime\prime}$ such that
$\pi_1^{\prime\prime} <_\st \pi_1$.  In order to conclude, it suffices
to show that any product
$\pi_1^{\prime\prime}\dots\pi_{s+1}^{\prime\prime}$ such that
$\pi_1^{\prime\prime} <_\st \pi_1$ is in $V_\st^{s+1}$.  In such a
product, $\pi_2^{\prime\prime}\dots\pi_{s+1}^{\prime\prime}$ may be
rewritten (using $V_\st^s=V^s$) as the sum of an element of $V_\st^{s-1}
\subseteq V_\st^{s}$ and a linear combination of terms
$\pi_2^{\prime\prime\prime}\dots\pi_{s+1}^{\prime\prime\prime}$ with
$\pi_2^{\prime\prime\prime}\le_\st\dots\le_\st\pi_{s+1}^{\prime\prime\prime}$. 
Thus, such a product
$\pi_1^{\prime\prime}\pi_2^{\prime\prime}\dots\pi_{s+1}^{\prime\prime}$
may be rewritten as the sum of an element of $V^s = V_\st^{s}$ and a
linear combination of terms
$\pi_1^{\prime\prime}\pi_2^{\prime\prime\prime}
\dots\pi_{s+1}^{\prime\prime\prime}$ with
$\pi_2^{\prime\prime\prime}\le_\st\dots\le_\st\pi_{s+1}^{\prime\prime\prime}$. 
But, since $\pi_1^{\prime\prime} <_\st \pi_1$, we have
$\rk\pi_1^{\prime\prime} < \rk\pi_1$ and, by the inductive hypothesis,
$\pi_1^{\prime\prime}\pi_2^{\prime\prime\prime}\dots
\pi_{s+1}^{\prime\prime\prime} \in V_\st^{s+1}$.  It follows that $\pi
\in V_\st^{s+1}$ as required.  This completes the proof.\qed

\begin{sublemma} \label{pre-st-fam-gen} -- In the notation above,
$V^s_\st = V^s$, for each $s\in {\mathbb N}$.  \end{sublemma}

\proof 
We proceed by induction on $s$, the cases $s=0,1,2$ being already
proved.  Assume that the result is true for some integer $s \ge 2$.  We
must show that $V^{s+1} \subseteq V^{s+1}_\st$; that is, we must show
that $\Pi^{s+1} \subseteq V^{s+1}_\st$.  Let $\pi$ be an element of
$\Pi^{s+1}$. If $\pi \in \Pi^{s} \subseteq V^s$,
we have $\pi \in V^s=V^s_\st \subseteq V^{s+1}_\st$,  by the inductive 
hypothesis; and so 
it suffices to consider the case where $\pi=\pi_1\dots\pi_{s+1}$, where 
$\pi_1,\dots,\pi_{s+1} \in \Pi$.  By the inductive hypothesis,
$\pi_2\dots\pi_{s+1} \in V_\st^s$ and so this product 
may be written as a sum of
$\omega \in V^{s-1}_\st$ and a linear combination of standard monomials
of length $s$.  It follows that we may rewrite $\pi$ as the sum of
$\pi_1\omega \in V^s=V^s_\st \subseteq V^{s+1}_\st$ and a linear
combination of elements of $\Pi^{s+1}$ of the form
$\pi_1\pi_2^\prime\dots\pi_{s+1}^\prime$, where
$\pi_2^\prime,\dots,\pi_{s+1}^\prime$ are elements of $\Pi$ such that
$\pi_2^\prime\le_\st\dots\le_\st\pi_{s+1}^\prime$.  By the 
inductive hypothesis and Lemma \ref{pre-pre-st-fam-gen}, such elements
are in $V^{s+1}_\st$.  This completes the proof.\qed

\begin{subproposition} \label{standard-mon-basis} -- 
Let $A$ be a quantum graded A.S.L. on the set $\Pi$. The set of standard 
monomials on $\Pi$ form a $\Bbbk$-basis of $A$.
\end{subproposition}

\proof The set of standard monomials on $\Pi$ is free over $\Bbbk$ 
by hypothesis; and so  
we need only 
prove that this set 
generates $A$ as a $\Bbbk$-vector space. However, since $\Pi$ 
generates $A$
as a $\Bbbk$-algebra, it is enough to show that a product of elements of 
$\Pi$ is a linear
combination of standard monomials. This is a consequence of Lemma
\ref{pre-st-fam-gen}. \qed

\begin{subproposition} \label{GK-dimension-ASL} -- Let $A$ be a quantum graded 
A.S.L. on the set $\Pi$; then $\GKdim A =\rk\Pi$.
\end{subproposition}

\proof We use the notation introduced before Lemma \ref{pre-pre-st-fam-gen}. 
In this notation, $V:=V^1$ is the span of $\Pi \cup \{1\}$ and
$V^s$ is the usual $s$-th power of the vector 
space $V$ for $s\geq 2$.
Since $1 \in V$ and since $V$ is a finite dimensional vector space 
generating $A$ as a $\Bbbk$-algebra, we have 
$\GKdim A = \overline{\lim}_{n \rightarrow\infty} \log_n 
(\dim V^n)$. Thus, by Lemma \ref{pre-st-fam-gen}, 
\[
\GKdim A = \overline{\lim}_{n \rightarrow\infty} \log_n  (\dim V^n_\st). 
\]
Let us denote by $M$ the number of subsets of $\Pi$ which are totaly ordered for
$<_\st$ and maximal (with respect to inclusion) for this property.
These sets all have rank less than or equal to that of $\Pi$ and, clearly, 
for $n \in {\mathbb N}$: 
\[
\dim V^n_\st \le M \left(\begin{array}{cc} n + \rk\Pi  \cr \rk\Pi
\end{array}\right),
\]
so that $\GKdim A \le \rk\Pi$.
On the other hand, there is a totally ordered
subset of $\Pi$ of rank $\rk\Pi$. Hence, 
\[
\dim V^n_\st \ge  \left(\begin{array}{cc} n + \rk\Pi  \cr \rk\Pi
\end{array}\right);
\] 
and so $\GKdim A \ge \rk\Pi$. \qed\\

To each element $\alpha \in \Pi$ we associate the ideal $I_\alpha$ of $A$ 
generated by
those $\gamma$ in $\Pi$ such that $\gamma <_\st \alpha$ (with the convention 
that an
ideal generated by the empty set is $\langle 0 \rangle$). 

\begin{sublemma} \label{normality-lemma} -- Let $A$ be a quantum graded A.S.L. 
on the set $\Pi$. Any element $\alpha$ in $\Pi$ is normal modulo the ideal 
$I_\alpha$. 
In particular, any minimal element in $\Pi$ is normal. In addition, 
if $\Pi$ has a unique
minimal element $\alpha$ for $<_\st$, then $\alpha$ is a regular normal element
of $A$. 
\end{sublemma}

\proof  Condition (5) in the definition of a quantum graded A.S.L. shows that, 
for all
$\beta$ in $\Pi$, there exists a scalar $c_{\alpha\beta} \in \Bbbk^\ast$ 
such that
$\alpha\beta-c_{\alpha\beta}\beta\alpha \in I_\alpha$. Denoting by 
$\pi_\alpha$ the
canonical projection of $A$ onto $A/I_\alpha$, this shows that, 
$\pi_\alpha(\alpha)$
commutes up to non zero scalar with all the elements $\pi(\beta)$, 
for $\beta \in \Pi$.
However, the elements $\pi(\beta)$, with 
$\beta \in \Pi$, form a set of algebra 
generators
of $A/I_\alpha$, so $\alpha$ is normal modulo the ideal $I_\alpha$. 
Now, assume  that
$\alpha$ is the unique minimal element of $\Pi$. In this case,
$\alpha \le_\st \beta$, 
for all
$\beta\in\Pi$. Thus, if $\omega$ is a standard monomial,
$\alpha\omega$ is still a standard monomial.  Since 
the set standard monomials is a $\Bbbk$-basis of $A$, 
by Proposition
\ref{standard-mon-basis}, it follows
that left multiplication by $\alpha$ in $A$ is an injective map.
However, as we already
mentioned above, condition (5) in the definition of a quantum graded 
A.S.L. shows
that $\alpha$ commutes up to non-zero scalar with 
any standard monomial.
Thus, the injectivity of right  multiplication by $\alpha$ in $A$ follows from 
the
injectivity of left multiplication. We have shown that $\alpha$ is a regular 
normal
element of $A$. \qed

\subsection{Properties of $\Pi$-ideals.} \label{subsect-pi-ideals}

Let $\Pi$ be a set with a partial order $<_\st$. If $<_\tot$ is a total order
on $\Pi$, we say that $<_\tot$ {\em respects} 
$<_\st$ if, for $\alpha,\beta\in\Pi$,
$\alpha <_\st \beta \Longrightarrow \alpha<_\tot\beta$. Notice that, for a 
finite poset
$(\Pi,<_\st)$, there always exists a total order on $\Pi$ which respects 
$<_\st$.\\

Let $\Pi$ be a set with a partial order $<_\st$. If $\Omega$ is a subset of 
$\Pi$, we say that $\Omega$ is a {\em $\Pi$-ideal} provided it satisfies 
the following condition: if $\alpha \in \Omega$ and if $\beta \in \Pi$,
with 
$\beta\le_\st\alpha$, then $\beta\in\Omega$. (Notice that what we call
a $\Pi$-ideal is called an {\it ideal in $\Pi$} in [BV]; see [BV; Proposition
5.1.].)

\begin{sublemma} \label{N-S-G} -- Let $A$ be a quantum graded A.S.L. on the 
poset $(\Pi,<_\st)$ and let $\Omega$ be a $\Pi$-ideal. In addition, 
let $<_\tot$
be any total
order on $\Pi$ which respects $<_\st$. Then the elements of $\Omega$, ordered 
by $<_\tot$ form a normalising sequence of generators of the ideal 
$\langle\Omega\rangle$ generated by $\Omega$ in $A$.
\end{sublemma}

\proof Let $\Omega=\{\omega_1 , \dots ,\omega_s\}$ with $\omega_1 <_\tot
\dots <_\tot \omega_s$. For $1 \le i \le s$, consider
$\alpha\in\Pi$ such that $\alpha <_\st
\omega_i$. Then, since $\Omega$ is a $\Pi$-ideal, $\alpha\in\Omega$. 
Moreover, since
$<_\tot$ respects $<_\st$, we also have $\alpha <_\tot \omega_i$. 
Thus, $\alpha \in
\{\omega_1,\dots,\omega_{i-1}\}$ (this set being empty if $i=1$). 
Thus, $I_{\omega_i}
\subseteq \langle\omega_1,\dots,\omega_{i-1}\rangle$. 
But then $\omega_i$ is normal modulo
$\langle\omega_1,\dots,\omega_{i-1}\rangle$ since, by 
Lemma \ref{normality-lemma}, it is normal
modulo $I_{\omega_i}$. \qed

\begin{sublemma} \label{ASL-noeth} -- A quantum graded A.S.L. is noetherian and
satisfies polynomial growth ((PG) for short) in the sense of [Lev; 5.4].
\end{sublemma}

\proof Let $A$ be a quantum graded A.S.L. on the poset $(\Pi,<_\st)$ and let
$<_\tot$ be any total order on $\Pi$ which respects $<_\st$. 
Since $\Pi$ is clearly a
$\Pi$-ideal, by Lemma \ref{N-S-G}, the elements of $\Pi$ ordered by $<_\tot$ 
form a nomalising sequence in $A$. Since $A/\langle\Pi\rangle \cong \Bbbk$ is 
noetherian, [ATV; Lemma 8.2]
shows that $A$ is noetherian. In the same way, since $A/\langle\Pi\rangle 
\cong \Bbbk$ 
satisfies (PG), [Lev; Proposition 5.6] shows that $A$ satisfies (PG).\qed\\

A nice feature of $\Pi$-ideals is the fact that they can be described in a very 
simple
way in terms of the standard monomials. This description is given in the next
proposition. If $\omega \in \Pi$, we say that a standard monomial $\mu\neq 1$ 
{\em involves} 
$\omega$ if $\mu=\omega_1 \dots \omega_s$, where 
$\omega_1,\dots,\omega_s \in \Pi$,
with $\omega_1 \le_\st \dots \le_\st \omega_s$, and $\omega \in
\{\omega_1,\dots,\omega_s\}$. \\

The following remark will be useful. 

\begin{subremark} \label{sub-pi-ideal} -- 
\rm Let $A$ be a quantum graded A.S.L.
on the poset $(\Pi,<_\st)$ and let $<_\tot$ be any total order on $\Pi$
which respects $<_\st$. 
Let $\Omega=\{\omega_1,\dots,\omega_s\}$, for some $s\geq 2$, 
be a $\Pi$-ideal where
$\omega_1,\dots,\omega_s$ are elements of $\Pi$ such that
$\omega_1<_\tot\dots<_\tot\omega_s$ and let
$\Omega'=\{\omega_1,\dots,\omega_{s-1}\}$. 
Consider $\pi \in \Pi$ such that there
exists $1 \le i \le s-1$ with $\pi \le_\st \omega_i$. Then, 
$\pi \in \Omega$ since
$\Omega$ is a $\Pi$-ideal. Now, suppose that $\pi=\omega_s$. Then we have 
$\omega_s=\pi \le_\st \omega_i$. It follows that $\omega_s=\pi \le_\tot
\omega_i<_\tot\omega_s$, a contradiction. Thus, $\pi \in \Omega'$. 
We have proved that $\Omega'$ is also a $\Pi$-ideal.
\end{subremark}

\begin{subproposition} \label{basis-for-ideals-gen-by-pi-ideals} -- Let $A$ be
a quantum graded A.S.L. on the set $\Pi$ and let $\Omega$ be a $\Pi$-ideal. 
The set of standard monomials involving an element of $\Omega$ form a 
$\Bbbk$-basis of $\langle\Omega\rangle$.
\end{subproposition}

\proof Let $<_\tot$ be any total order on $\Pi$ which respects $<_\st$.
Notice that the set of standard monomials involving an element of 
$\Omega$ is free by condition (3) of the definition of a quantum graded A.S.L. 
and that it is included in $\langle\Omega\rangle$. Thus it remains to prove
that this set generates $\langle\Omega\rangle$ as a $\Bbbk$-vector space.

We proceed by induction on the cardinality of $\Omega$.  
First, suppose that $\Omega$ is a $\Pi$-ideal with $|\Omega|=1$, say
$\Omega = \{\omega_1\}$. Then $\omega_1$ is a minimal element of $\Pi$
with respect to $\le_\st$, and $\langle\Omega\rangle=\omega_1 A$, by
Lemma \ref{N-S-G}. 
Let  $V$ be the subspace of $A$ generated by those standard monomials
involving $\omega_1$. Clearly,  $V \subseteq \langle\Omega\rangle$. 
Conversely, suppose that 
$x \in \langle\Omega\rangle$; and so there is an element $x_1
\in A$ such that $x=\omega_1 x_1$.  Set $x_1=c_1y_1+\dots + c_ty_t$,
where $c_1,\dots, c_t$ are scalars and $y_1, \dots, y_t$ standard
monomials; so that $x=\omega_1x_1=c_1\omega_1y_1+\dots +c_t\omega_1y_t$. 
Thus, it is enough to show that, if $\pi_1, \dots,\pi_s$ are elements in
$\Pi$ such that $\pi_1 \le_\st \dots \le_\st \pi_s$, then
$\omega_1\pi_1\dots\pi_s \in V$.  Of course, if $\omega_1 \le_\st
\pi_1$, this is trivial.  Thus, we assume that $\omega_1 \not\le_\st
\pi_1$.  Since $\omega_1$ is minimal in $\Pi$, it follows that
$\omega_1$ and $\pi_1$ are not comparable with respect to $\le_\st$. 
Thus, condition (4) of the definition of a quantum graded A.S.L. 
entails $\omega_1\pi_1=0$ and then $\omega_1\pi_1\dots\pi_s \in V$.  We
have shown that the result is true when $|\Omega|=1$. 

Next, suppose that  $s  > 1$, and assume that the result is
true for all the $\Pi$-ideals of cardinality less than $s$. Set 
$\Omega=\{\omega_1,\dots,\omega_s\}$ where $\omega_1,\dots,\omega_s$ are
elements of $\Pi$ such that $\omega_1<_\tot\dots<_\tot\omega_s$.  Then,
$\langle\Omega\rangle
=\omega_1A+\dots+\omega_sA$, by Lemma \ref{N-S-G}.  
Let $V$ be the $\Bbbk$-subspace of $A$
generated by those standard monomials involving an element of $\Omega$. 
Again, it is clear that $V \subseteq \langle\Omega\rangle$.  Now, let
$\Omega'=\{\omega_1,\dots,\omega_{s-1}\}$ and denote by $V'$ the
$\Bbbk$-subspace of $A$ spanned by those standard monomials involving an
element of $\Omega'$.  In fact, $\Omega'$ is a $\Pi$-ideal as we noticed
in Remark \ref{sub-pi-ideal}.  Thus, the inductive hypothesis gives
$\langle\Omega'\rangle =\omega_1A+\dots+\omega_{s-1}A=V'$.  Now, let $x
\in \langle\Omega\rangle$.  Then there exist $x_1,\dots,x_s \in A$ such
that $x=\omega_1x_1+\dots+\omega_sx_s$.  We want to show that $x \in V$. 
Since $\omega_1x_1+\dots+\omega_{s-1}x_{s-1} \in \langle\Omega'\rangle
=V' \subseteq V$, it is enough to show that $\omega_sx_s \in V$. 
Of course, since $x_s$ is a linear combination of standard monomials, it
is enough to show that, if $\pi_1, \dots,\pi_t$ are elements in $\Pi$
such that $\pi_1 \le_\st \dots \le_\st \pi_t$, then
$\omega_s\pi_1\dots\pi_t \in V$.  
First, assume that $\omega_s \not\ge_\st \pi_j$ for each $j$ with $1 \le
j \le t$. 
In this case, if $\omega_s$
and $\pi_1$ are comparable with respect to $\le_\st$ then  $\omega_s
<_\st \pi_1$; and so $\omega_s\pi_1\dots\pi_t$ is a standard monomial
which is 
obviously in $V$.  
Next, consider the case where 
$\omega_s$ and $\pi_1$ are not comparable
with respect to $\le_\st$. In this case, by condition (4) of the definition
of quantum graded A.S.L., $\omega_s\pi_1$ is a linear combination of
standard monomials involving an element $\lambda \in \Pi$ such that
$\lambda <_\st \omega_s$.  
Any such $\lambda$ is in $\Omega$, since $\Omega$ is a $\Pi$-ideal; and
so $\lambda \in \Omega'$, since $\lambda \neq \omega_s$.  It follows that
$\omega_s\pi_1 \in \langle\Omega'\rangle$ and then that
$\omega_s\pi_1\dots\pi_t \in \langle\Omega'\rangle \subseteq V$.  It
remains to consider the case where there exists an integer $i$, with 
$1 \le i
\le t$, such that $\omega_s \ge_\st \pi_i$.  Let $j$ denote the greatest 
such integer.  Thus, we have $\pi_1 \le_\st \dots \le_\st \pi_j \le_\st
\omega_s$.  By condition (5) of the definition of quantum graded
A.S.L., for $1 \le i \le j$, there exists a scalar $c_i \in \Bbbk^\ast$
such that $\omega_s\pi_i-c_i\pi_i\omega_s$ is a linear combination of
standard monomials involving an element $\lambda \in \Pi$ such that
$\lambda <_\st \omega_s$.  As we have already mentioned, this implies that
each such $\lambda$ is in $\Omega'$.  Thus, 
$\omega_s\pi_i-c_i\pi_i\omega_s \in \langle\Omega'\rangle$, for $1 \le i
\le j$.  Setting 
$c=c_1 \dots c_j$, it follows that 
$\omega_s\pi_1\dots\pi_j-c\pi_1\dots\pi_j\omega_s \in
\langle\Omega'\rangle$ and, as a consequence,
$\omega_s\pi_1\dots\pi_t-c\pi_1\dots\pi_j\omega_s\pi_{j+1}\dots\pi_t \in
\langle\Omega'\rangle$.  Since $\langle\Omega'\rangle=V' \subseteq V$,
it remains to prove that $\pi_1\dots\pi_j\omega_s\pi_{j+1}\dots\pi_t \in
V$.  Again, we isolate two cases.  First, assume that $\omega_s$ and
$\pi_{j+1}$ are comparable.  Then, by definition of $j$, we have
$\omega_s <_\st \pi_{j+1}$.  In this case,
$\pi_1\dots\pi_j\omega_s\pi_{j+1}\dots\pi_t$ is a standard monomial
that is obviously in $V$.  Next, assume that $\omega_s$ and
$\pi_{j+1}$ are not comparable.  Then, as we have already noticed,
$\omega_s\pi_{j+1}$ is a linear combination of standard monomials
involving an element $\lambda \in \Omega'$.  Hence,
$\pi_1\dots\pi_j\omega_s\pi_{j+1}\dots\pi_t \in  \langle\Omega'\rangle
\subseteq V$.  The proof is complete.\qed

\begin{subcorollary} \label{quotients-graded-ASL} -- Let $A$ be a quantum
graded A.S.L. on the set $\Pi$ and let $\Omega$ be a $\Pi$-ideal. 
Suppose that  $p \, : \, A
\longrightarrow A/\langle\Omega\rangle$ is the canonical projection.
Then the $\Bbbk$-algebra $A/\langle\Omega\rangle$ is a quantum graded 
A.S.L. on the set $p(\Pi\setminus\Omega)$.
\end{subcorollary}

\proof Since $\langle\Omega\rangle$ 
is generated by homogeneous elements, it is a graded 
ideal and $A/\langle\Omega\rangle$ 
inherits an ${\mathbb N}$-grading. By Proposition
\ref{basis-for-ideals-gen-by-pi-ideals} and the linear independence of standard
monomials, $\Pi \cap \langle\Omega\rangle 
= \Omega$. In fact $p(\Pi\setminus\Omega) \cong
\Pi\setminus\Omega$ and we may equip $p(\Pi\setminus\Omega)$ with the partial
order inherited from that of $\Pi$. Clearly, the elements of 
$p(\Pi\setminus\Omega)$
are homogeneous elements of positive degree generating $A/\langle\Omega\rangle$
as a $\Bbbk$-algebra.
In addition, the family of standard monomials on $p(\Pi\setminus\Omega)$ is free.
The fact that conditions (4) and (5) are satisfied for 
$A/\langle\Omega\rangle$ is then clear. 
\qed

\section{The AS-Cohen-Macaulay property.}

Throughout this section, $\Bbbk$ is a field.

\subsection{AS-Cohen-Macaulay algebras.} 

In this subsection, we recall the notions of AS-Cohen-Macaulay and AS-Gorenstein
algebras as defined in [JZ]. We also recall some results about these notions
that we will need latter. \\

Throughout this subsection, $A$ stands for a noetherian ${\mathbb
N}$-graded connected $\Bbbk$-algebra; that is, $A=\oplus_{i \ge 0} A_i$ is an
${\mathbb N}$-graded $\Bbbk$-algebra such that $A_0=\Bbbk$.  It follows that $A$
is a finitely generated $\Bbbk$-algebra and that it is 
locally finite in the sense that $A_n$ is finite dimensional as a 
$\Bbbk$-vector space for each $n \in {\mathbb N}$.  For any integer $n$, let
$A_{\ge n}=\oplus_{i \ge n} A_i$. We denote the unique maximal graded 
ideal, $A_{\ge 1}$, of $A$ by ${\mathfrak m}$. 
The category of graded left $A$-modules and homogeneous 
homomorphisms of degree
zero is denoted $\GrMod(A)$.  If $M=\oplus_{i \in {\mathbb Z}} M_i$ is
an object of $\GrMod(A)$, then for any integer $n$, the graded module
$M(n)$ is the $n$-th shift of $M$; that is, $M(n)=\oplus_{i \in {\mathbb
Z}} M(n)_i$ where $M(n)_i=M_{n+i}$ for $i \in {\mathbb Z}$. Let $M,N$
be objects in $\GrMod(A)$, then 
$\Hom_{\GrMod(A)}(M,N)$ will stand for the abelian
group of homogeneous homomorphisms of degree $0$. The $\Ext$ groups associated
to $\Hom_{\GrMod(A)}$ are then denoted $\Ext_{\GrMod(A)}^i$. In addition, we put
\[
\HOM_A(M,N):=
\displaystyle\bigoplus_{m \in {\mathbb Z}} \Hom_{\GrMod(A)}(M,N(m)).
\] 
The $\Ext$ groups associated to $\HOM_A$ are then denoted $\EXT_A^i$. Thus, if 
$M,N$ are objects in $\GrMod(A)$, 
\[
\EXT_A^i (M,N)
\cong \displaystyle\bigoplus_{m \in {\mathbb Z}} \Ext_{\GrMod(A)}^i(M,N(m)).
\]
We now define the local cohomology functors. If $M$ is an object in $\GrMod(A)$,
we define its {\em torsion submodule}, $\Gamma_{\mathfrak m}(M)$, by
\[
\Gamma_{\mathfrak m}(M):=\{m \in M \tq \exists\, n \in {\mathbb N},\;
A_{\ge n}m=0\}.
\] 
It is easy to check that $\Gamma_{\mathfrak m}$ (and restriction of morphisms)
defines a left exact functor. We denote the 
associated right derived functors by 
$H_{\mathfrak m}^i$. Then, if $M$ is an object in $\GrMod(A)$, 
\[  
H_{\mathfrak m}^i(M) \cong \lim_{\longrightarrow} \EXT^i_A (A/A_{\ge n},M).
\]  
We write $A^\circ$ for the opposite algebra of $A$. Of course, $A^\circ$ is 
again a
noetherian connected ${\mathbb N}$-graded algebra and for any integer $i$, the 
$i$-th
local cohomology functor associated to $A^\circ$ will be denoted $H_{\mathfrak
m^\circ}^i$.   

\begin{subdefinition} -- {\rm [JZ; Definition 0.1, Definition 0.2].} \\
The algebra $A$ is said to be 
{\em AS-Cohen-Macaulay} if there exists $n \in {\mathbb N}$
such that
\[
\forall\, i \in {\mathbb N}, \, i \neq n,  H_{\mathfrak m}^i(A) = H_{{\mathfrak
m}^\circ}^i(A) = 0.
\] 
The algebra $A$ is said to be {\em AS-Gorenstein} 
if it satisfies the following
conditions.\\ (1) We have $\injdim_A A = \injdim_{A^\circ} A = n < \infty$ (the
injective dimension being measured in $\GrMod(A)$ and $\GrMod(A^\circ)$
respectively).\\ (2) There exists an integer $\ell$ such that 
\[
\EXT^i_A (\Bbbk,A) \cong \EXT^i_{A^\circ} (\Bbbk,A) 
\cong \left\{ \begin{array}{cc} 0 & \mbox{ for } i \neq n, \cr \Bbbk(\ell) 
& \mbox{ for } i=n. \end{array} \right.
\]
\end{subdefinition}

\begin{subremark} -- \rm
Notice that, by [Lev; Lemma 3.3], in the definition of the 
AS-Gorenstein property, we might equally 
well use the usual injective dimension.
\end{subremark}

The following definition appeared in [Z; p. 392].

\begin{subdefinition} -- Let $B$ be a
noetherian connected ${\mathbb N}$-graded $\Bbbk$-algebra.  If for every
non-simple graded prime factor $B/P$ of $B$ there is a
non-zero homogeneous normal element in $(B/P)_{\ge 1}$, then we
say $B$ has {\em enough normal elements}.  
\end{subdefinition}

As we will see in Remark \ref{ref-co-case}, if $A$ has enough normal elements,
it behaves in a similar way to commutative noetherian local rings 
with respect to the AS-Gorenstein and AS-Cohen-Macaulay properties. 

\begin{subremark} -- \label{ASL-ENE} \rm \\
1. As mentioned in [Z; p.  392], if a noetherian
connected ${\mathbb N}$-graded $\Bbbk$-algebra $B$ has a normalising
sequence $x_1,\dots,x_n$ of elements in $B_{\ge 1}$ such that
$B/(x_1,\dots,x_n)$ is finite dimensional as a $\Bbbk$-vector space, then
$B$ as enough normal elements. \\
2. Clearly, a quantum graded A.S.L. is connected and it is noetherian by
Lemma \ref{ASL-noeth}. Hence,
from the previous point and Lemma \ref{N-S-G}, it follows at once that
any quantum graded A.S.L. has enough normal elements. 
\end{subremark}

The next result is due to Yekutieli and Zhang (see [YZ]). 
(Here, we take the convention that the Gelfand-Kirillov dimension of the
zero module is $-\infty$.) 

\begin{subtheorem} \label{borne-sup-hl} -- Assume that $A$ has
enough normal elements. 
Then, for each finitely generated graded left $A$-module
$M$,
\[
\sup\{i \in {\mathbb N} \tq H_{\mathfrak m}^i(M) \neq 0\}=\GKdim_A M < + \infty.
\]
\end{subtheorem}

\proof For unexplained terminology, the reader is referred to [YZ]. 
By the conjunction of [YZ; Theorem 5.13] and [YZ; Theorem 5.14], the
algebra $A$ is 
graded-Auslander (that is, $A$ has a graded-Auslander 
balanced dualizing complex $R$)  and for every 
finitely generated graded left $A$-module $M$,
we have $\Cdim M=\GKdim_A M$ (see [YZ; 4.8] and [YZ; 3.10]). 
It remains to apply [YZ; 4.14] to conclude that 
$\sup\{i \in {\mathbb N} \tq H^i_{\mathfrak m}(M) \neq 0\}
=\GKdim_A M < + \infty$, 
for every finitely generated 
graded left $A$-module $M$.  \qed\\

Given a finitely generated graded left $A$-module $M$, the previous theorem 
provides an
upper bound for the set of integers $i$ such that $H_{\mathfrak m}^i(M) \neq 0$. 
To obtain a lower bound, we need the notion of depth.

\begin{subdefinition} -- Let $M$ be a finitely generated graded left $A$-module, 
the depth of $M$ is defined by
\[
\depth_A M :=\inf\{i \in {\mathbb N} \tq \EXT_A^i(\Bbbk,M)\neq 0\} 
\in {\mathbb N} \cup \{+\infty\}.
\]
\end{subdefinition}

We list the results we need about depth in Lemma \ref{borne-inf-hl}.  We
also give a sketch of proof for those well-known results since we have
not been able to locate them in the literature. 

\begin{sublemma} \label{borne-inf-hl} -- \\
(i) Let $M$ be a finitely generated graded left $A$-module. Then 
\[
\depth_A M =\inf\{i \in {\mathbb N} \tq H_{\mathfrak m}^i(M) \neq 0\}.
\] 
(ii) If $0 \longrightarrow L \longrightarrow M \longrightarrow N \longrightarrow
0$ is a short exact sequence in $\GrMod(A)$ of finitely generated graded left 
$A$-modules, and if $\ell,m,n$ denote 
the depth of $L,M,N$, respectively, then 
\[
\ell \ge \min \{m,n+1\},\; m \ge \min\{\ell,n\},\; n \ge \min\{\ell -1,m\}.
\]
(iii) If $x \in A$ is a regular normal homogeneous element of positive degree,
then $\depth_{A/\langle x \rangle} A/\langle x \rangle=\depth_A A -1$.
\end{sublemma}

\proof Point (i) is proved as in the commutative case (see [W; Theorem 4.6.3]).\\
(ii) This follows immediately from examination 
of the long exact sequence of local
cohomology groups associated with the short exact sequence
$0 \longrightarrow L \longrightarrow M \longrightarrow N \longrightarrow
0$.\\
(iii) By the graded analogue of Rees' Lemma, we have 
\[ 
\forall \, i \in {\mathbb N},\qquad
\EXT_{A/\langle x \rangle}^i (\Bbbk,A/\langle x \rangle) 
\cong \EXT_A^{i+1} (\Bbbk,A).
\]
The result follows at once.\qed\\

\begin{subremark} \label{versus-grothendieck} -- \rm
Assume that $A$ has enough normal
elements and let $M$ be a non-zero finitely generated graded left $A$-module. 
By Theorem \ref{borne-sup-hl} and Lemma \ref{borne-inf-hl}, 
$\{i \in {\mathbb N} \tq H_{\mathfrak m}^i(M) \neq 0\}$ is non-empty,  
contains the integers $\depth_A M$ and $\GKdim_A M$ and any $i$ in this set 
satisfies
$\depth_A M \le i \le \GKdim_A M$. 
This result is similar to the classical vanishing theorem
of Grothendieck for commutative noetherian local rings 
(see [BH; Theorem 3.5.7]). 
\end{subremark}

\begin{subremark} \label{ref-co-case} -- \rm
The noetherian connected ${\mathbb N}$-graded algebras with enough
normal elements behave in a similar way to commutative noetherian
local rings from the point
of view of the AS-Cohen-Macaulay and the AS-Gorenstein properties.
Indeed, if $A$ has enough normal elements, then $A$ 
is AS-Gorenstein if and only if
$\injdim_A A = \injdim_{A^\circ} A < \infty$, 
see [Z; Proposition 2.3(2)],  and,
if $A$ is AS-Gorenstein, then it is AS-Cohen-Macaulay (use [Z; Theorem 3.1(1)],
[Z; Proposition 2.3(2)] and Remark \ref{versus-grothendieck}).
\end{subremark}

The following result is [YZ~; Lemma 4.5, p. 29], in conjunction with 
[JZ~; Lemma 5.8]. It will be useful later. 

\begin{subproposition} \label{change-of-ring-and-lc} -- 
Assume that $A$ has enough normal elements. 
Let $M$ be a graded left $A$-module. Let $I$ be 
a graded ideal in ${\rm Ann}_A
M$ and set ${\mathfrak n}={\mathfrak m}/I$. Then 
$H_{\mathfrak m}^i(M)$ and $H_{\mathfrak n}^i(M)$ are isomorphic as 
$A/I$-modules, for $i \in {\mathbb N}$. 
\end{subproposition}

\subsection{The AS-Cohen-Macaulay property for a Quantum graded A.S.L..}

The aim of this section is to prove that, if $A$ is a quantum graded
A.S.L.  on a so-called {\it wonderful} poset $\Pi$, then $A$ is
AS-Cohen-Macaulay.  \\

First, we prove a lemma which will be of central importance. 

\begin{subremark} \label{cns-ascm-gkdim-prof} \rm --  Let $A$ be a noetherian
${\mathbb N}$-graded connected $\Bbbk$-algebra and assume that $A$ has enough
normal elements. Let $K$ be an homogeneous ideal of $A$, with $K \neq A$.
Then the ring $A/K$ is AS-Cohen-Macaulay if and only if
\[
\depth_A A/K = \GKdim_A A/K = \GKdim_{A^\circ} A/K = \depth_{A^\circ} A/K.
\] 
(Note that the central equality is always true, by [KL; 5.4].) 
This remark is an easy consequence of Remark
\ref{versus-grothendieck} and Proposition \ref{change-of-ring-and-lc}.
\end{subremark}

\begin{sublemma} \label{central-lemma-CM} -- 
Let $A$ be a noetherian ${\mathbb N}$-graded connected $\Bbbk$-algebra
and assume that $A$ has enough normal elements.\\
(1) Assume that 
$A$ has polynomial growth ((PG) for short) in the sense of 
[Lev; 5.4].
Let $x \in A$ be a normal homogeneous element of positive degree which is not a 
zero divisor. 
Then, the ring $A/xA$ is AS-Cohen-Macaulay if and only if the ring $A$ 
is AS-Cohen-Macaulay.\\  (2) Suppose that $I$ and 
$J$ are homogeneous ideals of $A$ 
such that:\\ (i) $\GKdim_A A/I = \GKdim_A A/J =
\GKdim_A A$,\\ 
(ii) $\GKdim_A A/(I+J) = \GKdim_A A - 1$,\\ (iii) $I \cap J = (0)$,\\
(iv) the rings $A/I$ and $A/J$ are AS-Cohen-Macaulay.\\ Then, the ring $A$ is
AS-Cohen-Macaulay if and only if the ring $A/(I+J)$ is AS-Cohen-Macaulay.
\end{sublemma}

\proof (1) According to [Lev; 5.7], $\GKdim_A A/(x)=\GKdim_A A -1$ and
$\GKdim_{A^\circ} A/(x)=\GKdim_{A^\circ} A -1$. According to point (iii) in
Lemma \ref{borne-inf-hl} and Proposition \ref{change-of-ring-and-lc}, 
$\depth_A A/(x) = \depth_{A/(x)} A/(x) = \depth_A A -1$ and 
$\depth_{A^\circ} A/(x)=\depth_{A^\circ} A -1$. Thus, the
result follows at once from Remark \ref{cns-ascm-gkdim-prof}.\\ 
(2) Because of the hypothesis
(iii), there is a short exact sequence of left-$A$-modules (respectively
right-$A$-modules)
\begin{equation} \label{sec} 0 
\stackrel{}{\longrightarrow} A
\stackrel{f}{\longrightarrow} A/I \oplus A/J 
\stackrel{g}{\longrightarrow} A/(I+J)
\stackrel{}{\longrightarrow} 0, 
\end{equation} where, for $a,b \in A$, $f(a)=(a+I,-a+J)$ and $g((a+I,b+J))
=a+b+I+J$.
Indeed, clearly
$\ker f =I \cap J=(0)$, while 
$g$ is onto and $g\circ f=0$. Moreover, for $a,b \in A$,
if $0=g((a+I,b+J))=a+b+I+J$, then there exist $i\in I$ and $j \in J$ such that 
$a+b=i+j$.
Thus, $(a+I,b+J)=(x+I,-x+J)=f(x)$ 
where $x=a-i=j-b$. From these short exact sequences,
using point (ii) of Lemma \ref{borne-inf-hl}, we derive 
\begin{equation} \label{sec-prof-gauche}
\begin{array}{c}
\depth_A A \ge \min\{\depth_A A/I \oplus A/J,\depth_A A/(I+J)+1\},\cr
\depth_{A^\circ} A \ge \min\{\depth_{A^\circ} A/I \oplus A/J,\depth_{A^\circ}
A/(I+J)+1\},
\end{array}
\end{equation}
\begin{equation} \label{sec-prof-droit}
\begin{array}{c}
\depth_A A/(I+J) \ge \min\{\depth_A A -1,\depth_A A/I \oplus A/J\},\\
\depth_{A^\circ} A/(I+J) \ge \min\{\depth_{A^\circ} A -1,\depth_{A^\circ} A/I \oplus
A/J\}.
\end{array}
\end{equation} Using [KL; 5.4], hypotheses (i) and (ii) may be rewritten\\ (i)
$\GKdim_A A/I = \GKdim_A A/J = \GKdim_A A 
= \GKdim_{A^\circ} A = \GKdim_{A^\circ} A/I
= \GKdim_{A^\circ} A/J$,\\   (ii) $\GKdim_A A/(I+J) = \GKdim_A A - 1 =
\GKdim_{A^\circ} A -1 = \GKdim_{A^\circ} A/(I+J)$.\\ Because of Remark
\ref{cns-ascm-gkdim-prof}, hypothesis (iv) may be rewritten\\ (iv)
$\depth_A A/K = \GKdim_A A/K = \GKdim_{A^\circ} A/K = \depth_{A^\circ} A/K$, where
$K=I$ or $J$.\\  But, by [KL; 5.1], $\GKdim_A (A/I \oplus A/J)=\max\{\GKdim_A A/I,
\GKdim_A A/J\}$. Thus, as a consequence of hypothesis (i),
$\GKdim_A (A/I \oplus A/J)=\GKdim_A A$. In addition, point (ii) of Lemma 
\ref{borne-inf-hl} shows that $\depth_A
(A/I \oplus A/J)\ge\min\{\depth_A A/I, \depth_A A/J\}$. Thus, hypotheses (i) and (iv)
give $\depth_A (A/I \oplus A/J)\ge\GKdim_A A$. Since, by 
\ref{versus-grothendieck} 
the depth of a non-zero finitely generated graded module
is bounded by its GK-dimension, it follows that 
$\depth_A (A/I \oplus A/J)=\GKdim_A (A/I \oplus A/J)=\GKdim_A A$. Clearly, the same
holds with $A^\circ$ instead of $A$. Thus, we have proved
\begin{equation} \label{prof-oplus}
\begin{array}{c}
\depth_A (A/I \oplus A/J)=\GKdim_A (A/I \oplus A/J)=\GKdim_A A \cr
\depth_{A^\circ} (A/I \oplus A/J)=\GKdim_{A^\circ} (A/I \oplus A/J)=\GKdim_{A^\circ} A.
\end{array}
\end{equation} Now, assume that $A/(I+J)$ is AS-Cohen-Macaulay. By
\ref{cns-ascm-gkdim-prof}, this means that  $\depth_A A/(I+J) = \GKdim_A A/(I+J) =
\GKdim_{A^\circ} A/(I+J) =
\depth_{A^\circ} A/(I+J)$. Thus, using hypothesis (ii), this is equivalent to 
$\depth_A A/(I+J) = \GKdim_A A -1 = \GKdim_{A^\circ} A -1 = \depth_{A^\circ} A/(I+J)$.
Thus, from relations
(\ref{sec-prof-gauche}) and (\ref{prof-oplus}), it follows that 
$\depth_A A \ge \GKdim_A A$ and $\depth_{A^\circ} A \ge \GKdim_{A^\circ} A$. But, as
we noticed before, the depth is bounded above 
by the GK-dimension, thus $\depth_A A =
\GKdim_A A = \GKdim_{A^\circ} A = \depth_{A^\circ} A$. This means  
that the ring $A$ is AS-Cohen-Macaulay, by Remark 
\ref{cns-ascm-gkdim-prof} .

A similar argument, using (\ref{sec-prof-droit}) instead of
(\ref{sec-prof-gauche}) shows that, if the ring $A$ is
AS-Cohen-Macaulay, then the ring $A/(I+J)$ is AS-Cohen-Macaulay.  \qed\\

Now, we need to introduce combinatorial tools.  For this we essentially
follow [BV; Chapter 5.D].  An important notion, in what follows, is that of
a {\em wonderful poset}.  For the definition of such posets, the reader
is referred to [BV; Chapter 5.D]. Note that any distributive lattice is a 
wonderful poset (see [BV; p. 58]).  We will also use Lemma 5.13 of [BV] which
contains all the properties of wonderful posets that we need.  Finally, let
$\Pi$ be a set partially ordered by $<_\st$; following [BV; Chap. 5.A], 
to any subset $\Sigma$ of
$\Pi$ we associate two $\Pi$-ideals as follows.  The {\em $\Pi$-ideal
generated by $\Sigma$} is the smallest $\Pi$-ideal containing $\Sigma$; 
that is, the set $\{\xi\in\Pi \tq
\exists\,\sigma\in\Sigma,\,\xi\le_\st\sigma\}$.  The {\em $\Pi$-ideal
cogenerated by $\Sigma$} is the greatest $\Pi$-ideal disjoint from
$\Sigma$; that is, the set $\{\xi\in\Pi \tq
\forall\, \sigma\in\Sigma,\,\xi\not\ge_\st\sigma\}$. 

\begin{subtheorem} \label{critirion-ASL-ASCM} -- If $A$ is a quantum
graded A.S.L. on a wonderful poset $(\Pi,<_\st)$, 
then $A$ is AS-Cohen-Macaulay.  \end{subtheorem}

\proof Recall that, by Remark \ref{ASL-ENE}, 
a quantum graded A.S.L. is a noetherian ${\mathbb N}$-graded connected 
$\Bbbk$-algebra with enough normal elements.
The proof is by induction on $|\Pi|$.  The result is true if
$|\Pi|=1$ since, in this case, $A$ is a (commutative) polynomial ring in
one indeterminate, which is well known to be AS-Cohen-Macaulay.  Now, we
suppose the result is true up to cardinality $n \ge 1$ and we consider a
quantum graded A.S.L.  $A$ on a wonderful poset $\Pi$ of cardinality
$n+1$. Two cases may occur according to the number of minimal elements in $\Pi$. 

First, assume that $(\Pi,<_\st)$ has a single minimal element $\xi$.  By
Lemma \ref{normality-lemma}, $\xi$ is a regular normal element of $A$. 
In addition, Lemma \ref{ASL-noeth} shows that $A$ satisfies (PG). 
Moreover, Corollary \ref{quotients-graded-ASL}, [BV; Lemma 5.13 (b)] and
the inductive hypothesis show that, $A/\xi A$ is AS-Cohen-Macaulay.  So,
by part (1) of Lemma \ref{central-lemma-CM}, $A$ is AS-Cohen-Macaulay. 

Now, assume that the minimal elements of $\Pi$ are $\xi_1,\dots,\xi_k$,
with $k \ge 2$.  Let $\Omega$ be the $\Pi$-ideal cogenerated by $\xi_1$ and
$\Psi$ be the $\Pi$-ideal cogenerated by $\xi_2,\dots,\xi_k$.  We let
$I=\langle\Omega\rangle$ 
and $J=\langle\Psi\rangle$.  Consider $x \in I \cap J$ and assume that 
$x$
is non-zero.  Then, by Proposition \ref{standard-mon-basis} we may write 
$x=\lambda_1x_1+\dots +
\lambda_sx_s$, with $s \ge 1$, where $\lambda_1,\dots,\lambda_s$ are non-zero
scalars and $x_1,\dots,x_s$ are standard monomials.  Using Proposition
\ref{basis-for-ideals-gen-by-pi-ideals}, we see that for each $i \in
\{1,\dots,s\}$, the term 
$x_i$ must involve an element $\omega_i \in \Omega$ and
an element $\psi_i \in \Psi$.  But, $x_i$ being standard, $\omega_i$
and $\psi_i$ are comparable for $<_\st$ and, $\Omega$ and $\Psi$ being
$\Pi$-ideals, the smaller element of $\omega_i$ and $\psi_i$ must
be in $\Omega \cap \Psi$.  However, [BV; 5.13 (e)(iii)] shows that $\Omega
\cap \Psi=\emptyset$; a contradiction.  Thus, we have 
$I \cap J=\langle 0 \rangle$. 
We have $I+J=\langle \Omega \cup \Psi \rangle$ and clearly $\Omega \cup \Psi$
is a $\Pi$-ideal. 
By [BV; 5.13(e)], $\Pi\setminus\Omega$, $\Pi\setminus\Psi$ and 
$\Pi\setminus(\Omega\cup\Psi)$ are wonderful posets and we have 
$\rk(\Pi\setminus\Omega)=\rk(\Pi\setminus\Psi)=\rk\Pi$ and 
$\rk(\Pi\setminus(\Omega\cup\Psi))=\rk\Pi -1$. Now, from 
Corollary \ref{quotients-graded-ASL} and Proposition \ref{GK-dimension-ASL}, 
we deduce that $\GKdim_A A/I=\GKdim_A A/J=\GKdim_A A$ 
and $\GKdim_A A/(I+J)=\GKdim_A A -1$ (see [KL; 5.1(c)]). In addition,
the induction hypothesis gives that the rings $A/I$, $A/J$ and $A/(I+J)$ are
AS-Cohen-Macaulay. So, we are in position to apply Lemma \ref{central-lemma-CM}
which shows that the ring $A$ is AS-Cohen-Macaulay.  
The proof is complete.  \qed

\section{Examples of quantum graded A.S.L.}

The aim of this section is to show that quantum grassmannians and quantum
determinantal rings are
quantum graded A.S.L. \\

For our purposes, we need to define quantum grassmannians in a 
slightly more general setting than usual; that is, over an arbitrary 
commutative domain. We also need to introduce quantum analogues of 
coordinate rings of Schubert varieties. This is done in the first subsection.
In the second subsection, 
we exhibit useful bases for quantum grassmannians, called
{\it standard bases}, in this general setting, thus extending results from
[KLR] where such bases were constructed in the case where the base ring
is a field. Using this material, we show in Subsections 3 and 4 that quantum
grassmannians over a field are examples of quantum graded A.S.L. Finally, in
Subsection 5, we deduce that quantum determinantal rings over a field are 
quantum graded A.S.L.\\

Throughout the rest of this work, we fix two positive integers $m,n$. 

\subsection{Quantum grassmannians.} \label{subsec-q-grass}

Recall the definition of the quantum analogue, ${\cal O}_q(M_{m,n}(\Bbbk))$, 
of the coordinate ring of $m \times n$ matrices with entries in 
the field $\Bbbk$.

\begin{subdefinition} -- \label{def-q-mat}
Let $\Bbbk$ be a field and $q$ be a non-zero element of $\Bbbk$.
We denote by ${\cal O}_q(M_{m,n}(\Bbbk))$ the $\Bbbk$-algebra generated
by $mn$ indeterminates $X_{ij}$, for 
$1\le i \le m$ and $1\le j \le n$, 
subject to the following relations:
\[
\begin{array}{ll} 
X_{ij}X_{il}=qX_{il}X_{ij}, 
& \mbox{ for }1\le i \le m, \mbox{ and } 1\le j<l \le n\: ; \\ 
X_{ij}X_{kj}=qX_{kj}X_{ij}, 
& \mbox{ for }1\le i<k \le m, \mbox{ and } 1\le j \le n \: ; \\ 
X_{ij}X_{kl}=X_{kl}X_{ij}, 
& \mbox{ for } 1\le k<i \le m, \mbox{ and } 1\le j<l \le n \: ; \\ 
X_{ij}X_{kl}-X_{kl}X_{ij}=(q-q^{-1})X_{il}X_{kj}, 
& \mbox{ for } 1\le i<k \le m, \mbox{ and } 1\le j<l \le n. 
\end{array}
\]
\end{subdefinition}

More generally, we make the following definition. 

\begin{subdefinition} -- Let $A$ be a commutative integral domain and denote
by $\Bbbk$ its field of fractions. For any unit $u \in A$, we define 
${\cal O}_u(M_{m,n}(A))$ to be the $A$-subalgebra of 
${\cal O}_u(M_{m,n}(\Bbbk))$ generated over $A$ by $\{X_{ij}\}$, 
for $1 \le i \le m$ and
$1 \le j \le n$.
\end{subdefinition}

\begin{subremark} \label{q-grass-over-ring-as-ioe} -- \rm 
Let $A$ be a commutative domain whose field of fractions is denoted by
$\Bbbk$ and let $u$ be a unit of $A$.\\
1. It is well-known that ${\cal O}_u(M_{m,n}(\Bbbk))$ is an iterated skew
polynomial extension of $\Bbbk$ where the indeterminates are added in the order
$X_{11},\dots,X_{1n},X_{21},\dots,X_{2n},\dots,X_{m1},\dots,X_{mn}$.
As a well-known consequence, ${\cal O}_u(M_{m,n}(\Bbbk))$ is a noetherian
domain.\\
2. It is not difficult to see, as a consequence of the previous point,
that ${\cal O}_u(M_{m,n}(A))$ is an iterated Ore extension of the ring $A$
where, at each step, the endomorphism of the skew-derivation that is used is
in fact an automorphism.\\
3. It is well known that there is a $\Bbbk$-algebra isomorphism ${\rm Tr} \,
: \, {\cal O}_u(M_{m,n}(\Bbbk)) \longrightarrow {\cal O}_u(M_{n,m}(\Bbbk))$
such that  $X_{ij} \mapsto X_{ji}$, which we call the 
{\em transpose isomorphism}.
It follows at once that there is a similar transpose isomorphism of 
$A$-algebras ${\rm Tr}\, : \, {\cal O}_u(M_{m,n}(A)) \longrightarrow
{\cal O}_u(M_{n,m}(A))$. 
\end{subremark}

It follows from point 3 in Remark \ref{q-grass-over-ring-as-ioe} that, in what
follows, we may assume that $m \le n$ without loss of generality.  We do
this for the rest of the section.\\

If $A$ is a commutative domain and $u$ a unit in $A$, we associate to 
${\cal O}_u(M_{m,n}(A))$ the $m \times n$ matrix $X=(X_{ij})$ which we call the
{\em generic matrix} of ${\cal O}_u(M_{m,n}(A))$. \\

Let $1 \le t \le m$.
A pair $(I, J)$ is an {\em index pair (of cardinality $t$)} provided
that $I \subseteq \{1,\dots,m\}$ and $J\subseteq \{1,\dots,n\}$ with
$|I| = |J| = t$. 
If $I=\{i_1,\dots,i_t\}$ and $J=\{j_1,\dots,j_t\}$
with $1 \le i_1<\dots <i_t\le m$ and $1 \le j_1<\dots < j_t\le n$, we will
write $I=\{i_1<\dots <i_t\}$ and $J=\{j_1<\dots < j_t\}$. 
The set of all index pairs of cardinality $t$, where $t$ takes all possible
values in $\{1,\dots,m\}$ is denoted
$\Delta_{m,n}$ (or $\Delta$ if there is no possible confusion). To any
index pair $(I,J)$ of cardinality $t$ with $I=\{i_1<\dots <i_t\}$ and 
$J=\{j_1<\dots < j_t\}$ we associate the  {\em quantum $t \times t$
minor} 
$[I|J]$ of $X$. This is the element of ${\cal O}_u(M_{m,n}(A))$ defined by
\[
[I|J]:=\sum_{\sigma \in {\mathfrak S}_t} (-u)^{\ell(\sigma)} 
X_{i_1,j_{\sigma(1)}} \dots X_{i_t,j_{\sigma(t)}}
\]
(here, ${\mathfrak S}_t$ is the symmetric group on $\{1,\dots,t\}$ and $\ell$
is the usual length function on ${\mathfrak S}_t$).\\

In the same way, an  {\em index set} $J$ is 
a subset $J$ of $m$ pairwise distinct 
elements of $\{1,\dots,n\}$. If $J=\{j_1,\dots,j_m\}$
with $1 \le j_1 < \dots < j_m \le n$, we will denote $J=\{j_1<\dots < j_m\}$.
The set of all index sets is denoted
$\Pi_{m,n}$ (or $\Pi$ if there is no possible confusion). \\
To any $J=\{j_1< \dots <j_m\} \in \Pi$, we associate the maximal quantum minor
$[J]:=[K|J]$ of $X$, where $K=\{1 <\dots < m\}$, hence
\[
[J]=\sum_{\sigma \in S_m} (-u)^{\ell(\sigma)} X_{1,j_{\sigma(1)}} \dots
X_{m,j_{\sigma(m)}}.
\]
The {\em standard partial order}, $\le_\st$, on $\Pi$ is 
defined as follows:
\[
\{i_1<\dots<i_m\} \le_\st \{j_1<\dots<j_m\} \qquad\mbox{iff}\qquad i_s \le j_s
\qquad\mbox{for each}\qquad 1 \le s \le m.
\]
Clearly, $\Pi$ can be identified with the subset of maximal
quantum minors of ${\cal  O}_u(M_{m,n}(A))$. Thus, in what follows, we
also denote by $\Pi$ the set of maximal 
quantum minors of ${\cal  O}_u(M_{m,n}(A))$.

\begin{subdefinition} -- \label{def-q-grass}
Let $A$ be a commutative domain and $u$ be a unit 
in $A$. The {\em quantum grassmannian} ${\cal O}_u(G_{m,n}(A))$ is 
the $A$-subalgebra of 
${\cal O}_u(M_{m,n}(A))$ generated by $\{[I],\, I \in \Pi\}$. 
\end{subdefinition}

The following remark will be of constant use in what follows. We will often
use it without explicit reference.

\begin{subremark} \label{homo-grass} -- \rm 
Let $A$ and $B$ be commutative domains, let 
$u \in A$ and $v \in B$ be units
and let $f \, : \, A \longrightarrow B$ 
be a morphism of rings such that $f(u)=v$.
By using Remark \ref{q-grass-over-ring-as-ioe} and the universal property
of Ore extensions, it is easy to see that there is a unique morphism of
rings $g \, :\, {\cal O}_u(M_{m,n}(A)) \longrightarrow {\cal O}_v(M_{m,n}(B))$
such that $g(X_{ij})=X_{ij}$ for $1 \le i \le m$ and $1 \le j \le n$ which
extends $f$.
Clearly, for $I \in \Pi$, we then have $g([I])=[I]$. Hence, $g$ induces by
restriction a morphism of rings $h \, :\, {\cal O}_u(G_{m,n}(A))
\longrightarrow {\cal O}_v(G_{m,n}(B))$ such that $h([I])=[I]$ for any
$I\in\Pi$. It is easy to see that, if $f$ is injective (respectively, 
surjective), then both $g$ and $h$ are injective (respectively, surjective).
\end{subremark}

\begin{subremark} \label{various-q-matrices} -- \rm 
In what follows, we will be mainly interested in the following cases.\\
1. The case where $A={\mathbb Z}[t^{\pm 1}]$ and $u=t$. This case will
be, in fact, of central importance, thus we fix the notation ${\mathcal A}
={\mathbb Z}[t^{\pm 1}]$ for the rest of this work.\\
2. The case where $A=\Bbbk$ is a field and $q$ is any non-zero element in 
$\Bbbk$. Then, ${\cal O}_q(M_{m,n}(\Bbbk))$ 
is the usual quantum deformation of the coordinate ring on the 
variety $M_{m,n}(\Bbbk)$. Moreover, ${\cal O}_q(G_{m,n}(\Bbbk))$ 
is the usual quantum deformation of the 
homogeneous coordinate ring of the grassmannian $G_{m,n}(\Bbbk)$.
\end{subremark}

In order to prove that quantum grassmannians have a straightening law,
we will need to use certain quantum analogues
of {\it Schubert varieties}. These are quotients of quantum grassmannians
and the key to obtaining straightening relations in quantum grassmannians is 
to find good bases for these quantum Schubert varieties.

\begin{subdefinition} -- Let $A$ be a commutative domain and $u$ be a unit in 
$A$. To any $\gamma \in \Pi$ we associate the set 
$\Omega_\gamma:=\{\pi \in \Pi \tq
\pi \not\ge_\st \gamma\}$ and we set
\[
I_{\gamma,A}:=\langle\Omega_\gamma\rangle  \qquad\mbox{and}\qquad 
{\cal O}_u(G_{m,n}(A))_\gamma:={\cal O}_u(G_{m,n}(A))/I_{\gamma,A}.
\]
We denote by $\Theta_{\gamma,A}$ the canonical projection 
\[
\Theta_{\gamma,A} \, : \, {\cal O}_u(G_{m,n}(A)) \longrightarrow {\cal
O}_u(G_{m,n}(A))_\gamma.
\]
The quotient ring ${\cal O}_u(G_{m,n}(A))_\gamma$ will be called the 
{\em quantum Schubert variety (over $A$) associated with $\gamma$}.
\end{subdefinition}

\subsection{Standard bases for quantum grassmannians.}

\begin{subdefinition} -- 
Let $A$ be any commutative domain and $u$ be any invertible
element of $A$. A {\em standard monomial} on $\Pi$ is an element of 
${\cal O}_u(G_{m,n}(A))$ which is either  equal to $1$ or is of the form 
$[I_1]\dots[I_s]$, with $I_1,\dots,I_s\in\Pi$ and 
$I_1 \le_\st\dots\le_\st I_s$.
\end{subdefinition}

\begin{subremark} \label{basis-grass-over-field} -- \rm
Recall from [KLR; Corollary 2.8] that, if $A$ is any field and $u$ is any
nonzero element of $A$ then 
the standard monomials on $\Pi$ form an $A$-basis of 
${\cal O}_u(G_{m,n}(A))$.
\end{subremark}

We want to extend Remark \ref{basis-grass-over-field} to any 
commutative integral domain $A$. 

\begin{sublemma} -- \label{gen-fam-over-integers-Laurent}
The set of standard monomials in ${\cal O}_t(G_{m,n}({\mathbb Z}[t^{\pm 1}]))$ 
is a ${\mathbb Z}[t^{\pm 1}]$-generating set of 
${\cal O}_t(G_{m,n}({\mathbb Z}[t^{\pm 1}]))$.
\end{sublemma}

\proof Let $x$ be a non zero element of ${\cal O}_t(G_{m,n}({\mathbb
Z}[t^{\pm 1}]))$.  We consider the natural embeddings ${\mathbb
Z}[t^{\pm 1}] \hookrightarrow {\mathbb Q}(t)$ and ${\mathbb Z}[t^{\pm
1}] \hookrightarrow {\mathbb C}[t^{\pm 1}]$ and the associated injective
morphisms of rings $h \,:\, {\cal O}_t(G_{m,n}({\mathbb Z}[t^{\pm 1}]))
\hookrightarrow {\cal O}_t(G_{m,n}({\mathbb Q}(t)))$ and $h' \,:\, {\cal
O}_t(G_{m,n}({\mathbb Z}[t^{\pm 1}])) \hookrightarrow {\cal
O}_t(G_{m,n}({\mathbb C}[t^{\pm 1}]))$. 

Remark
\ref{basis-grass-over-field} applied with $A={\mathbb Q}(t)$ and $u=t$,
shows that $h(x)$ may be written as a linear combination 
of standard monomials in ${\cal O}_t(G_{m,n}({\mathbb Q}(t)))$ with
coefficients from  ${\mathbb Q}(t)$.  
It follows at once that there are elements 
$k(t),k_1(t),\dots,k_s(t) \in {\mathbb Z}[t^{\pm 1}]
\setminus\{0\}$ and standard monomials $M_1,\dots,M_s$ in ${\cal
O}_t(G_{m,n}({\mathbb Z}[t^{\pm 1}]))$ such that $k(t)x =\sum_{i=1}^s
k_i(t) M_i \in {\cal O}_t(G_{m,n}({\mathbb Z}[t^{\pm 1}]))$.  Of course,
we may assume that $k(t),k_1(t),\dots,k_s(t)$ have no common irreducible
factor in ${\mathbb Z}[t^{\pm 1}]$ and hence in ${\mathbb C}[t^{\pm
1}]$.  By applying $h'$ to this relation we obtain a relation
$k(t)h'(x)=\sum_{i=1}^s k_i(t) M_i \in {\cal O}_t(G_{m,n}({\mathbb C}
[t^{\pm 1}])$.  Now, consider $u \in {\mathbb C}^\ast$: to the natural
map ${\it ev}_u \, : \, {\mathbb C}[t^{\pm 1}] \longrightarrow {\mathbb
C}$, $t \mapsto u$ we may associate an epimorphism ${\cal
O}_t(G_{m,n}({\mathbb C} [t^{\pm 1}])) \longrightarrow {\cal
O}_u(G_{m,n}({\mathbb C}))$.  Applying this map with $u$ a non-zero root
of $k(t)$ gives a nontrivial relation of linear dependence among
standard monomials in ${\cal O}_u(G_{m,n}({\mathbb C}))$ which violates
Remark \ref{basis-grass-over-field}.  Hence, there exists $c,d \in
{\mathbb Z}$ such that $k(t)=ct^d$.  As a consequence, we have the
relation $cx =\sum_{i=1}^s t^{-d}k_i(t) M_i$ in ${\cal
O}_t(G_{m,n}({\mathbb Z}[t^{\pm 1}]))$.  Here again, we may assume that
in ${\mathbb Z}$  
the greatest common divisor of $c$ and the non zero coefficients of the
$k_i(t)$, for $1 \le i \le s$, is one.  Now, assume that $p$ is a
prime divisor of $c$ in ${\mathbb Z}$.  The natural map ${\mathbb
Z}[t^{\pm 1}] \longrightarrow {\mathbb F}_p[t^{\pm 1}] \hookrightarrow
{\mathbb F}_p(t)$ gives rise to a map $h'' \, : \, {\cal
O}_t(G_{m,n}({\mathbb Z}[t^{\pm 1}])) \longrightarrow {\cal
O}_t(G_{m,n}({\mathbb F}_p(t)))$.  Applying $h''$ to the above relation
gives $0=\sum_{i=1}^s t^{-d}(k_i(t)+p{\mathbb Z}[t^{\pm 1}]) M_i$ in
${\cal O}_t(G_{m,n}({\mathbb F}_p(t)))$ where the right hand side is a
non trivial linear combination of standard monomials in ${\cal
O}_t(G_{m,n}({\mathbb F}_p(t)))$.  Again, this violates Remark
\ref{basis-grass-over-field}.  Thus, $c=\pm 1$ and we get the relation
$x =\sum_{i=1}^s (\pm 1) t^{-d}k_i(t) M_i$ in ${\cal
O}_t(G_{m,n}({\mathbb Z}[t^{\pm 1}]))$, as required.  \qed

\begin{subproposition} \label{set-gen-over-ring} --
Let  $A$ be  a commutative domain and let 
$u$ be an invertible element of $A$.
The set of standard monomials in ${\cal O}_u(G_{m,n}(A))$ is a $A$-generating
set of ${\cal O}_u(G_{m,n}(A))$.
\end{subproposition}

\proof Since $u$ is a unit in $A$, there is a natural ring homomorphism
${\mathbb Z}[t^{\pm 1}] \longrightarrow A$ such that $1 \mapsto 1_A$ and 
$t \mapsto u$. Let $f \, : \, {\cal O}_t(G_{m,n}({\mathbb Z}[t^{\pm 1}])) 
\longrightarrow {\cal O}_u(G_{m,n}(A))$ denote the induced morphism. 
Clearly, the set of products of the form $[I_1]\dots[I_s]$, 
with $I_1,\dots,I_s \in \Pi$ together with $1$ form a generating set for the 
$A$-module ${\cal O}_u(G_{m,n}(A))$. Hence it is enough to prove that any
such product is a linear combination with coefficients in $A$ of standard
monomials. Let $I_1,\dots,I_s \in \Pi$; we have $[I_1]\dots[I_s]
=f([I_1]\dots[I_s])$. By Lemma \ref{gen-fam-over-integers-Laurent}, in 
${\cal O}_t(G_{m,n}({\mathbb Z}[t^{\pm 1}]))$, we have an identity 
$[I_1]\dots[I_s]= \sum_{i=1}^d k_i(t)M_i$, where 
$k_i(t) \in {\mathbb Z}[t^{\pm 1}]$ and $M_i$ is a standard monomial, for
$1 \le i \le d$. Applying $f$ to this equality shows that in
${\cal O}_u(G_{m,n}(A))$, we have an identity
$[I_1]\dots[I_s]= \sum_{i=1}^d f(k_i(t))M_i$, where 
$f(k_i(t)) \in A$ and $M_i$ is a standard monomial, for $1 \le i \le d$.\qed

\begin{subproposition} \label{basis-over-ring} -- 
Let 
$A$ be a commutative domain and let $u$ be an invertible element of $A$. 
The set of standard monomials in ${\cal O}_u(G_{m,n}(A))$ is a $A$-basis of
${\cal O}_u(G_{m,n}(A))$.
\end{subproposition}

\proof Let $F$ be the field of fractions of $A$.
The natural embedding $A \hookrightarrow F$ induces an embedding
$\epsilon_A\, : \, {\cal O}_u(G_{m,n}(A)) \hookrightarrow 
{\cal O}_u(G_{m,n}(F))$ 
with $\epsilon_A([I_1]\dots[I_s])=[I_1]\dots[I_s]$, 
for $I_1 \le_\st \dots \le_\st I_s$ in 
$\Pi$. The linear independence 
of standard monomials in 
${\cal O}_u(G_{m,n}(A))$ thus follows from the linear independence of standard
monomials in ${\cal O}_u(G_{m,n}(F))$ (see Remark 
\ref{basis-grass-over-field}). On the other hand, 
the standard monomials generate ${\cal O}_u(G_{m,n}(A))$ 
as an $A$-module, by Proposition
\ref{set-gen-over-ring}. \qed

\subsection{Straightening relations in quantum grassmanians.}

The aim of this subsection is to establish that quantum grassmannians
have a straightening law. In order to complete this task, we will make use
of the quantum Schubert varieties introduced earlier.

\begin{subremark} \label{omega-gamma-pi-ideal} --  \rm 
Let $A$ be a commutative domain and $u$ a unit of $A$.\\
1. We say that a standard monomial $[I_1]\dots[I_s]$ in 
${\cal O}_u(G_{m,n}(A))$
{\em involves} an element $[J]$, with $J \in \Pi$, 
if $J\in\{I_1,\dots,I_s\}$.\\
2. It is clear from the definition that 
the subset 
$\Omega_\gamma$ is a $\Pi$-ideal for any
$\gamma\in\Pi$. As a consequence, a standard
monomial $[I_1]\dots[I_s]$ in ${\cal O}_u(G_{m,n}(A))$ involves an 
element of $\Omega_\gamma$ if an only if $I_1 \in \Omega_\gamma$.
\end{subremark}

Recall that ${\mathcal A}={\mathbb Z}[t^{\pm 1}]$.
Our main aim in this section is to exhibit useful ${\mathcal A}$-bases 
for quantum Schubert varieties over ${\mathcal A}$  
when $u=t$. In order to complete this task, we need some commutation
relations established in the works [HL(1)] and [HL(2)] by Hodges and 
Levasseur. Since the point of view and the notation of these papers are quite 
different from those of this work, we recall some facts from them, for the 
convenience of the reader. Let $r \in {\mathbb N}^\ast$. Recall that, 
for any field $\Bbbk$ and any non-zero element $q \in \Bbbk$, the {\em
quantum special linear group} is defined by  
${\mathcal O}_q(SL_r(\Bbbk)):=
{\mathcal O}_q(M_r(\Bbbk))/\langle D_q -1\rangle$,
where ${\mathcal O}_q(M_r(\Bbbk))={\mathcal O}_q(M_{r,r}(\Bbbk))$ and $D_q$
is the quantum determinant of ${\mathcal O}_q(M_r(\Bbbk))$; that is, 
its unique $r \times r$ quantum minor.\\

The setting of [HL(1)] and [HL(2)] is the following.
Fix an integer $r\geq 2$ and a complex number $q$ which is not a root of 
unity. Let ${\mathfrak g}$ denote the Lie algebra 
${\mathfrak s}{\mathfrak l}_r({\mathbb C})$ and $G$ the algebraic group 
$SL_r({\mathbb C})$. Let  $U_q({\mathfrak g})$ be the quantum 
universal enveloping algebra of ${\mathfrak g}$ over ${\mathbb C}$ as defined 
in [HL(1)~; \S 1.2] and let ${\mathbb C}_q[G]$ be  
the restricted dual associated with $U_q({\mathfrak g})$. 
According to [HL(1); Theorem 1.4.1], the algebra 
${\mathbb C}_q[G]$ is nothing but 
${\mathcal O}_{q^2}(SL_r({\mathbb C}))$ in our notation. 
Now, for any integer $k$, with $1 \le k \le r-1$, and any element
$w$ in the Weyl group $W \cong {\mathfrak S}_r$ associated with 
${\mathfrak s}{\mathfrak l}_r({\mathbb C})$,
we consider the matrix coefficient $c_{k,w}^+ \in {\mathbb C}_q[G]$ defined
in [HL(1); 1.5]. As noted in [HL(1); 1.5], and making no distinction between
a minor of ${\mathcal O}_{q^2}(M_r({\mathbb C}))$ and its image in 
${\mathcal O}_{q^2}(SL_r({\mathbb C}))$, we have
\[
c_{k,w}^+ = [I|K] \in {\mathbb C}_q[G]={\mathcal O}_{q^2}(SL_r({\mathbb C})), 
\] 
where $K=\{1,\dots,k\}$ and $I$ is the index-set obtained from $\{w(1),\dots,
w(k)\}$ after reordering. Then, for $x,y \in W$,
applying [HL(2); Proposition 1.1] with $1 \le i=j=k
\le r-1$ and $t = {\rm id} \in {\mathfrak S}_r$, we get an identity
\begin{equation} \label{HL-cr-relation}
c_{k,y}^+c_{k,x}^+ - (q^2)^e c_{k,x}^+c_{k,y}^+ 
= \zzum{u \in W}{x<_k ux}{uy <_k y} g_u(q) c_{k,ux}^+ c_{k,uy}^+ , 
\end{equation}
where $e \in {\mathbb Z}$, and $g_u(q) \in {\mathbb C}$, for $u \in W$ 
such that $x<_k ux$ and $uy <_k y$.
(For the definition of the relation $<_k$ in $W$, see [HL(2)].) 

Now, set 
$I=\{i_1< \dots <i_k\}$ and $J=\{j_1< \dots <j_k\}$
with $I,J \subseteq \{1,\dots,r\}$ and $k\leq r-1$. Let
$K=\{1,\dots,k\}$. Once translated
in to our notation for minors, the relation (\ref{HL-cr-relation}) 
yields a relation
\begin{equation} \label{prefinal-cr}
[I|K][J|K]-(q^2)^e [J|K][I|K] = \zum{J<_\st L}{M <_\st I} g_{L,M} [L|K][M|K],
\end{equation}
in ${\mathcal O}_{q^2}(SL_r({\mathbb C})$, where $g_{L,M} \in {\mathbb C}$, for
$L,M \subseteq \{1,\dots,r\}$ such that $|L|=|M|=k$ 
while  $J<_\st L$ and $M <_\st I$. 
(Here again, we use the same notation for a minor in 
${\mathcal O}_{q^2}(M_r({\mathbb C}))$ and its image in 
${\mathcal O}_{q^2}(SL_r({\mathbb C}))$.)
On the other hand, by using the transpose automorphism of 
${\mathcal O}_{q^2}(M_r({\mathbb C}))$ (see [PW; 3.7.1(1)] and [PW; 4.1.1])
we can interchange rows and colums in 
(\ref{prefinal-cr}). Finally, it is well-known that, for $1 \le k < r$, 
the natural map
${\mathcal O}_{q^2}(G_{k,r}({\mathbb C})) \hookrightarrow 
{\mathcal O}_{q^2}(M_r({\mathbb C})) \longrightarrow
{\mathcal O}_{q^2}(SL_r({\mathbb C}))$ is an embedding. 
Hence, we end up with the
following important proposition.

\begin{subproposition} \label{HL-cr} -- We assume that 
$A={\mathbb C}$ and $u=q \in
{\mathbb C}$ is not a root of unity.
Let $I,J \in \Pi$. Then there exists a relation
\[
[I][J]-q^{e_{I,J}} [J][I] = \zum{J<_\st L}{M <_\st I} g_{L,M} [L][M],
\]
in ${\mathcal O}_q(G_{m,n}({\mathbb C}))$, where $e_{I,J} \in {\mathbb Z}$ and 
$g_{L,M} \in {\mathbb C}$, for
$L,M \subseteq \{1,\dots,n\}$ such that $|L|=|M|=m$ while 
$J<_\st L$ and $M <_\st I$.
\end{subproposition}

\begin{subtheorem} \label{theo-launois} --  
Let $A={\mathbb C}$ and suppose that 
$u=q \in {\mathbb C}$ is not a root of unity.
Let $\gamma \in \Pi$. Then the ideal $I_{\gamma,{\mathbb C}}$ of 
${\mathcal O}_q(G_{m,n}({\mathbb C}))$ is generated
by $\{\pi \in \Pi \tq \pi \not\ge_\st \gamma\}$ as a right ideal.  
\end{subtheorem}

\proof The ideal $I_{\gamma,{\mathbb C}}$ is generated, as a two-sided
ideal, by $\Omega_\gamma=\{\pi \in \Pi \tq \pi \not\ge_\st
\gamma\}$. Now, consider any total order on $\Pi$ that respects the
partial order $\le_\st$ on $\Pi$. 
Then the commutation relations in Proposition 
\ref{HL-cr} together with the fact that $\Omega_\gamma$ is a $\Pi$-ideal
show that the elements of $\Omega_\gamma$ ordered by this
total order form a normalising sequence in 
${\mathcal O}_q(G_{m,n}({\mathbb C}))$. The result follows. \qed\\

We will also need the following adaptation of [GL; Corollary 1.8].

\begin{subproposition} \label{rewrite-in-grass} -- 
Let  $A=K$ be a field and let $q$ be a nonzero element of
$K$. Let $I_1,\dots,I_s$ be elements of $\Pi$. In ${\mathcal O}_q(G_{m,n}(K))$,
the product $[I_1]\dots[I_s]$ can be expressed as a linear combination 
of standard monomials of the form $[J_1]\dots[J_s]$ such that 
$J_1 \le_\st I_1$.
\end{subproposition}  

\proof In this proof, we will use notation from [GL]. Recall from 
Remark \ref{basis-grass-over-field} that the standard monomials on $\Pi$ form a 
$K$-basis of ${\cal O}_q(G_{m,n}(K))$. Hence, the product $[I_1]\dots[I_s]$
is a linear combination of standard monomials in ${\cal O}_q(G_{m,n}(K))$
which must have the form $[J_1]\dots[J_s]$ with 
$J_1 \le_\st \dots \le_\st J_s \in \Pi$ 
as one may easily deduce using the 
standard grading of ${\cal O}_q(M_{m,n}(K))$. Now, consider the embedding 
${\cal O}_q(M_{m,n}(K)) \hookrightarrow {\cal O}_q(M_n(K))$ defined by 
$X_{ij} \mapsto X_{ij}$. It induces an embedding $ \iota \, : \,
{\cal O}_q(G_{m,n}(K)) 
\hookrightarrow {\cal O}_q(M_n(K))$ such that $[I] \mapsto [L|I]$ where $I$ is 
any element of $\Pi$ and $L=\{1,\dots,m\}$. Now, it is clear that, if 
$I,J \in \Pi$,
$I \le_\st J$ implies that $(L,I) \le (L,J)$ where $\le$ stands 
for the order $(\le_r,\le_c)$ 
defined in [GL; \S 1.2]. Thus, applying $\iota$ to the above
expression of $[I_1]\dots[I_s]$ as a linear combination of standard monomials
$[J_1]\dots[J_s]$ leads to an expression of $[L|I_1]\dots[L|I_s]$ in 
${\cal O}_q(M_n)$ as a linear combination of products $[L|J_1]\dots[L|J_s]$ 
whose associated bitableau 
\[
\left(
\begin{array}{cc}
L &  J_1 \cr \vdots & \vdots \cr L &  J_s
\end{array}
\right)
\]
is preferred (in the sense of [GL; \S 1.3]).
But, by [GL; Corollary 1.10], such products form a basis of ${\cal O}_q(M_n)$, 
hence we can compare this expression of $[L|I_1]\dots[L|I_s]$ with that given
by [GL; Corollary 1.8] to conclude that $(L,J_1) \le (L,I_1)$ which leads to 
$J_1 \le_\st I_1$. \qed

\begin{subcorollary} \label{I-gamma-as-vs} -- 
Let $A={\mathbb C}$ and suppose that $u=q \in {\mathbb C}$ is
not a root of unity.
Let $\gamma \in \Pi$.  The ideal $I_{\gamma,{\mathbb C}}$ 
of ${\mathcal O}_q(G_{m,n}({\mathbb C}))$ is the ${\mathbb C}$-span of 
the standard monomials which involve an element of $\Omega_\gamma$.   
\end{subcorollary} 

\proof We let $V_{\gamma,{\mathbb C}}$ denote the ${\mathbb C}$-span in
${\mathcal O}_q(G_{m,n}({\mathbb C}))$
of those standard monomials which involve an element of $\Omega_\gamma$.
Notice that, by Remark \ref{omega-gamma-pi-ideal},
a standard monomial involves an
element of $\Omega_\gamma$ if and only if it starts with an element of
$\Omega_\gamma$; that is, it has the form $[I_1]\dots[I_s]$ where  
$I_1,\dots,I_s\in\Pi$ with  $I_1\le_\st\dots\le_\st I_s$
and $I_1 \in \Omega_\gamma$. Clearly, $V_{\gamma,{\mathbb C}} \subseteq 
I_{\gamma,{\mathbb C}}$.

By Theorem \ref{theo-launois} and Remark
\ref{basis-grass-over-field}, any element in $I_{\gamma,{\mathbb C}}$ is a
linear combination of terms of the form $[I]\mu$ where 
$I\in\Omega_\gamma$ and $\mu$ is a standard monomial. Hence, to show that
$V_{\gamma,{\mathbb C}} \supseteq I_{\gamma,{\mathbb C}}$, it
is enough to show that any product $[I]\mu$ where 
$I\in\Omega_\gamma$ and $\mu$ is a standard monomial
lies in $V_{\gamma,{\mathbb C}}$. The case $\mu=1$ is trivial; and so 
we just have to consider products $[I][I_1]\dots[I_s]$ with 
$I_1,\dots,I_s\in\Pi$ and $I_1\le_\st\dots\le_\st I_s$.
By using Proposition \ref{rewrite-in-grass}, we see that  
$[I][I_1]\dots[I_s]$ is a linear combination of standard
monomials $[J_1]\dots[J_{s+1}]$, where $J_1
\le_\st I$. However,  $\Omega_\gamma$ is a $\Pi$-ideal; and so from 
$J_1 \le_\st I$ we deduce that $J_1 \in
\Omega_\gamma$. Hence $I_{\gamma,{\mathbb C}} \subseteq 
V_{\gamma,{\mathbb C}}$. \qed\\ 

We now investigate the case of quantum Schubert varieties over the
base ring ${\mathcal A={\mathbb Z}[t^{\pm 1}]}$. Let $q$ be an element
of ${\mathbb C}$ that is transcendental over ${\mathbb Q}$. There is a
canonical embedding 
${\mathbb Z}[t^{\pm 1}] \hookrightarrow {\mathbb C}$ sending
$t$ to $q$. This embedding induces an embedding 
\[
\varepsilon \, : \, {\cal O}_t(G_{m,n}({\mathbb Z}[t^{\pm 1}])) 
\hookrightarrow 
{\cal O}_q(G_{m,n}({\mathbb C})),
\]
such that, for $I \in \Pi$, $\varepsilon([I])=[I]$.
We denote by $I_{\gamma,{\mathcal A}}$ and
$I_{\gamma,{\mathbb C}}$ the ideals of 
${\cal O}_t(G_{m,n}({\mathbb Z}[t^{\pm 1}]))$
and  ${\cal O}_q(G_{m,n}({\mathbb C}))$, respectively, generated by the
elements of the $\Pi$-ideal $\Omega_\gamma$. In addition, we let 
$V_{\gamma,{\mathcal A}}$
denote the ${\mathbb Z}[t^{\pm 1}]$-span of those standard monomials in
${\cal O}_t(G_{m,n}({\mathbb Z}[t^{\pm 1}]))$ involving an element of 
$\Omega_\gamma$. (Recall from the proof of Corollary \ref{I-gamma-as-vs} that 
$V_{\gamma,{\mathbb C}}$
is the ${\mathbb C}$-span of those standard monomials in
${\cal O}_q(G_{m,n}({\mathbb C}))$ involving an element of $\Omega_\gamma$.)

\begin{sublemma} \label{about-I-gamma-A}-- In the notation above,
$I_{\gamma,{\mathcal A}}=\varepsilon^{-1}(I_{\gamma,{\mathbb C}})$
and $I_{\gamma,{\mathcal A}}=V_{\gamma,{\mathcal A}}$.
\end{sublemma}

\proof We have $V_{\gamma,{\mathcal A}} \subseteq I_{\gamma,{\mathcal A}}
\subseteq \varepsilon^{-1}(I_{\gamma,{\mathbb C}})
= \varepsilon^{-1}(V_{\gamma,{\mathbb C}})
= V_{\gamma,{\mathcal A}}$.
Indeed, the first inclusion is obvious, the second is easy since $\varepsilon$ 
is a ring homomorphism which sends $[I]$ to $[I]$ for any $I \in \Pi$, 
the first equality is 
Corollary \ref{I-gamma-as-vs} and the second equality is an easy consequence 
of the fact that the ${\mathbb Z}[t^{\pm 1}]$-basis of standard monomials in
${\cal O}_t(G_{m,n}({\mathbb Z}[t^{\pm 1}])$ 
is sent to the ${\mathbb C}$-basis 
of standard monomials in ${\cal O}_q(G_{m,n}({\mathbb C}))$ in the obvious
way and of the fact (due to the transcendance of $q$ over ${\mathbb Q}$) 
that a nonzero element in
${\mathbb Z}[t^{\pm 1}]$ is sent by $\varepsilon$ to a nonzero complex
number. Hence, we have proved the desired equalities.
\qed

\begin{subremark} \label{basis-schubert-C-Z-t} -- \rm
(i) It follows at once from Corollary \ref{I-gamma-as-vs} that, 
if $q \in {\mathbb C}$ is a non-root of unity, then 
the set of cosets of those standard monomials in 
${\cal O}_q(G_{m,n}({\mathbb C}))$
of the form $[I_1]\dots[I_s]$ such that $\gamma\le_\st I_1$ together with $1$
form a ${\mathbb C}$-basis of ${\cal O}_q(G_{m,n}({\mathbb C}))_\gamma$;\\
(ii) it then follows from  Lemma
\ref{about-I-gamma-A} that the set of cosets of those standard monomials in
${\cal O}_t(G_{m,n}({\mathbb Z}[t^{\pm 1}]))$ of the form 
$[I_1]\dots[I_s]$ such that $\gamma\le_\st I_1$ together with $1$
form a ${\mathbb Z}[t^{\pm 1}]$-basis of 
${\cal O}_t(G_{m,n}({\mathbb Z}[t^{\pm 1}]))_\gamma$. 
\end{subremark}

We now prove the existence of straightening relations in 
${\cal O}_q(G_{m,n}(A))$ for any commutative domain $A$.\\

Recall that, 
for any $\gamma\in\Pi$, we denote by 
$\Theta_{\gamma,{\mathbb Z}[t^{\pm 1}]}$ 
the canonical projection
\[
\begin{array}{ccc}
{\cal O}_t(G_{m,n}({\mathbb Z}[t^{\pm 1}])) 
& \stackrel{\Theta_{\gamma,{\mathbb Z}[t^{\pm 1}]}}{\longrightarrow} &
{\cal O}_t(G_{m,n}({\mathbb Z}[t^{\pm 1}]))_\gamma \cr
\end{array}.
\]

\begin{subtheorem} \label{straightening-relations-in-Z-Laurent} -- 
Let $I,J$ be elements of $\Pi$ which are not comparable with respect to 
$\le_\st$. Then, in ${\cal O}_t(G_{m,n}({\mathbb Z}[t^{\pm 1}]))$, 
the product $[I][J]$ may be rewritten as
\[
[I][J]=\sum_{i=1}^s k_i[I_i][J_i],
\]
with $k_i \in {\mathbb Z}[t^{\pm 1}]$, while 
$I_i,J_i \in \Pi$ with  $I_i <_\st I,J$ and $I_i \le_\st J_i$.
\end{subtheorem}

\proof From Proposition \ref{basis-over-ring}, we know that standard monomials 
form a ${\mathbb Z}[t^{\pm 1}]$-basis in 
the ring ${\cal O}_t(G_{m,n}({\mathbb Z}[t^{\pm 1}]))$. 
Thus, the product $[I][J]$ in 
${\cal O}_t(G_{m,n}({\mathbb Z}[t^{\pm 1}]))$ may be written as a linear 
combination of standard monomials. 
In addition, ${\cal O}_t(G_{m,n}({\mathbb Z}[t^{\pm 1}]))$ is a graded 
${\mathbb Z}[t^{\pm 1}]$-algebra in which the degree of the canonical 
generators is one. It follows that we have a relation
\[ 
[I][J]=\sum_{i=1}^s k_i[I_i][J_i],
\]
with $k_i \in 
{\mathbb Z}[t^{\pm 1}] \setminus \{0\}$ while $I_i,J_i \in \Pi$ and 
$I_i \le_\st J_i$. (Notice that 
the product $[I][J]$ is non zero, since ${\cal O}_t(G_{m,n}({\mathbb
Z}[t^{\pm 1}]))$ is a domain; and so the right hand side of the above
relation contains at least one term.)

Consider a maximal element $\gamma$ of $\{I_1,\dots,I_s\} \subseteq \Pi$ with 
respect to $\le_\st$. By applying $\Theta_{\gamma,{\mathbb Z}[t^{\pm 1}]}$ 
to the above relation, we obtain a relation
\[   
\Theta_{\gamma,{\mathbb Z}[t^{\pm 1}]}([I])
\Theta_{\gamma,{\mathbb Z}[t^{\pm 1}]}([J])
=\sum_{i=1}^s k_i\Theta_{\gamma,{\mathbb Z}[t^{\pm 1}]}([I_i])
\Theta_{\gamma,{\mathbb Z}[t^{\pm 1}]}([J_i]).
\] 
In the right hand side of this last relation, all the terms
$\Theta_{\gamma,{\mathbb Z}[t^{\pm 1}]}([I_i])
\Theta_{\gamma,{\mathbb Z}[t^{\pm 1}]}([J_i])$ such that 
$I_i \not\ge_\st \gamma$ are zero while the others form a non-empty family 
(it contains $\Theta_{\gamma,{\mathbb Z}[t^{\pm 1}]}([\gamma])
\Theta_{\gamma,{\mathbb Z}[t^{\pm 1}]}([J_i])$ for some $i$) which is free by 
Remark \ref{basis-schubert-C-Z-t}.
Thus, $\Theta_{\gamma,{\mathbb Z}[t^{\pm 1}]}([I])
\Theta_{\gamma,{\mathbb Z}[t^{\pm 1}]}([J]) \neq 0$ and we deduce that 
$I,J \ge_\st \gamma$. 

Thus, we have shown that, for any maximal element $\gamma$ of 
$\{I_1,\dots,I_s\}$, we have $I,J \ge_\st \gamma$. This entails that
$I_1,\dots,I_s \le_\st I,J$. Now, if $I_i=I$ or $I_i=J$ for some $i \in 
\{1,\dots,s\}$, since $I_i \le_\st I,J$, it follows that 
$I$ and $J$ must be comparable, which
is not the case. The proof is complete. \qed

\begin{subtheorem} \label{straightening-relations} -- Let $A$ be a commutative
domain, and let $u$ be a unit in $A$. Suppose that 
$I,J$ are elements of $\Pi$ which are not comparable with respect to $\le_\st$.
Then, in ${\cal O}_u(G_{m,n}(A))$, the product $[I][J]$ 
may be rewritten as
\[
[I][J]=\sum_{i=1}^s k_i[I_i][J_i],
\]
with 
$k_i \in A$ while 
$I_i,J_i \in \Pi$ with  $I_i <_\st I,J$ and $I_i \le_\st J_i$.
\end{subtheorem}

\proof There is a morphism of rings ${\mathbb Z} \longrightarrow A$ such that
$1 \mapsto 1_A$. Since $u$ is a unit in $A$, it induces a morphism of rings 
${\mathbb Z}[t^{\pm 1}] \longrightarrow A$ such that $t \mapsto u$. From this,
we get a morphism of rings 
${\mathcal O}_t(G_{m,n}({\mathbb Z}[t^{\pm 1}]))
\longrightarrow {\mathcal O}_u(G_{m,n}(A))$
sending $[I]$ to $[I]$. At this point, it is clear that
the desired relations in ${\mathcal O}_u(G_{m,n}(A))$
follow from the corresponding ones in 
${\mathcal O}_t(G_{m,n}({\mathbb Z}[t^{\pm 1}]))$ given by Theorem 
\ref{straightening-relations-in-Z-Laurent}. \qed

\subsection{Quantum grassmannians are quantum graded A.S.L.}

In order to show that quantum grassmannians are quantum graded A.S.L.,
it still remains to obtain good commutation relations for
them. This is the first goal of this section. 

\begin{subproposition} \label{ASL-pre-cr-over-complexes} -- 
Let $A={\mathbb C}$ and suppose that 
$u=q \in {\mathbb C}$ is not a root of unity.
Let $I,J \in \Pi$. Then there exists a relation
\[
[I][J]-q^{f_{I,J}} [J][I] = \sum k_{L,M} [L][M],
\]
in ${\mathcal O}_q(G_{m,n}({\mathbb C}))$, 
where the sum is taken over pairs $(L,M)$ of elements of $\Pi$
such that $L <_\st I$ and where $f_{I,J} \in {\mathbb Z}$ and 
$k_{L,M} \in {\mathbb C}$.
\end{subproposition}

\proof It is well-known that there exists an anti-isomorphism of
algebras $\alpha \, : \, {\mathcal O}_q(M_n(\mathbb C)) \longrightarrow  
{\mathcal O}_{q^{-1}}(M_n(\mathbb C))$ such that $X_{ij} \mapsto X_{ij}$, 
for
$1 \le i,j \le n$ (see, for example, [PW; 3.7]). In addition, for 
$I=\{i_1<\dots<i_t\}$ and $J=\{j_1<\dots<j_t\}$, 
with $1 \le i_1<\dots<i_t \le n$ and $1 \le j_1<\dots<j_t \le n$,
\[
\begin{array}{rcl}
\alpha([I|J]) &=& 
\alpha\left(\sum_{\sigma \in {\mathfrak S}_t} 
(-q)^{\ell(\sigma)}X_{i_1,j_{\sigma(1)}} \dots
X_{i_t,j_{\sigma(t)}}\right)
\cr &=& 
\sum_{\sigma \in {\mathfrak S}_t} 
(-q)^{\ell(\sigma)}X_{i_t,j_{\sigma(t)}}\dots X_{i_1,j_{\sigma(1)}}
\cr &=&  
\sum_{\sigma \in {\mathfrak S}_t} 
(-q)^{\ell(\sigma)}
X_{i_{\omega(1)},j_{\sigma\omega(1)}}\dots
X_{i_{\omega(t)},j_{\sigma\omega(t)}},
\end{array}
\]
where $\omega$ is the longest element of the symmetric group
${\mathfrak S}_t$, that is the element which sends $r$ to $t+1-r$, for
$1 \le r \le t$. As is well known, 
$\ell(\sigma)=\ell(\omega)-\ell(\sigma\omega)$, 
for all $\sigma \in {\mathfrak S}_t$.
Hence, we get
\[
\begin{array}{rcl}
\alpha([I|J])
& = & 
\sum_{\sigma \in {\mathfrak S}_t} 
(-q)^{\ell(\omega)-\ell(\sigma\omega)}
X_{i_{\omega(1)},j_{\sigma\omega(1)}}\dots 
X_{i_{\omega(t)},j_{\sigma\omega(t)}} \cr
& = & 
\sum_{\sigma \in {\mathfrak S}_t} 
(-q)^{\ell(\omega)-\ell(\sigma)}
X_{i_{\omega(1)},j_{\sigma(1)}}\dots X_{i_{\omega(t)},j_{\sigma(t)}} \cr
& = &
[I|J].
\end{array}
\]
(For the last equality, see [PW; 4.1.1].) It follows that $\alpha$
induces an anti-isomorphism 
$\beta \, : \, {\mathcal  O}_q(G_{m,n}(\mathbb C)) \longrightarrow 
{\mathcal  O}_{q^{-1}}(G_{m,n}(\mathbb C))$ such that, for any $I \in
\Pi$, $\beta([I])=[I]$. Now, let $I,J \in \Pi$. By Proposition \ref{HL-cr}
we have a relation 
\[
[I][J]-q^{e_{I,J}} [J][I] = \zum{J<_\st L}{M <_\st I} g_{L,M} [L][M],
\]
in ${\mathcal O}_{q^{-1}}(G_{m,n}({\mathbb C}))$, where 
$e_{I,J} \in {\mathbb Z}$ and 
$g_{L,M} \in {\mathbb C}$, for
$L,M \subseteq \{1,\dots,n\}$ such that $|L|=|M|=m$ and $J<_\st L, M <_\st I$.
Applying $\beta^{-1}$ to the identity gives us a new identity
\begin{equation}
[I][J] - q^{-e_{I,J}}[J][I] = \zum{J<_\st L}{M <_\st I} g_{L,M}^\prime [M][L].
\end{equation}
This completes the proof.\qed  

\begin{subproposition} \label{ASL-semipre-cr-over-complexes} --
Let $A={\mathbb C}$ and suppose that $u=q \in {\mathbb C}$
is not a root of unity.
Let $I,J \in \Pi$. Then there exists a relation
\[
[I][J]-q^{f_{I,J}} [J][I] = \sum k_{L,M} [L][M],
\]
in ${\mathcal O}_q(G_{m,n}({\mathbb C}))$,
where the sum is taken over pairs $(L,M)$ of elements of $\Pi$
such that $L <_\st I,J$ and either $L \le_\st M$ or $M <_\st L$ 
and where $f_{I,J} \in {\mathbb Z}$ and $k_{L,M} \in {\mathbb C}$.
\end{subproposition}

\proof Fix a pair $(I,J) \in \Pi\times\Pi$.
If $I$ and $J$ are not comparable for $\le_\st$,
then applying Theorem \ref{straightening-relations} to both $[I][J]$ and
$[J][I]$ shows that the result holds with, for example, $f_{I,J}=0$.
Now, suppose $I$ and $J$ are comparable. Of course, we may assume without
loss of generality that $I \le_\st J$. Now, Proposition
\ref{ASL-pre-cr-over-complexes} gives us a relation
$[I][J]-q^{f_{I,J}} [J][I] = \sum k_{L,M} [L][M]$,
where the sum is taken over pairs $(L,M)$ of elements of $\Pi$
such that $L <_\st I$ (and so $L <_\st J$ also holds, since $I \le_\st
J$) and where $f_{I,J} \in {\mathbb Z}$ and
$k_{L,M} \in {\mathbb C}$. In addition, by applying Theorem
\ref{straightening-relations} to each pair $(L,M)$ of incomparable
elements on the right hand side of this relation, 
we see that we can produce such a relation with the additional information
that, for each pair $(L,M)$ appearing, either $L \le_\st M$ or $M <_\st L$.
\qed

\begin{subtheorem} \label{ASL-cr-over-complexes} -- 
Let  $A={\mathbb C}$ and suppose that $u=q \in {\mathbb C}$
is not a root of unity.
Let $I,J \in \Pi$. Then there exists a relation
\[
[I][J]-q^{f_{I,J}} [J][I] = \sum k_{L,M} [L][M],
\]
in ${\mathcal O}_q(G_{m,n}({\mathbb C}))$, 
where the sum is taken over pairs $(L,M)$ of elements of $\Pi$
such that $L <_\st I,J$ and $L \le_\st M$ and where $f_{I,J} \in {\mathbb Z}$ 
and $k_{L,M} \in {\mathbb C}$.
\end{subtheorem}

\proof Let us introduce some notation.  We let $S$  be
the subset of those pairs $(L,M)$ in $\Pi\times\Pi$ such that $L \le_\st
M$  and let $T$ the subset of those pairs
$(L,M)$ in $\Pi\times\Pi$ such that 
$M <_\st L$.  In addition, for each pair $(I,J) \in
\Pi\times\Pi$, we let $E_{I,J}$ be the subset of those pairs $(L,M)$ in
$\Pi\times\Pi$ such that $L <_\st I,J$.  In this notation, Proposition
\ref{ASL-semipre-cr-over-complexes} shows that, for each pair $(I,J) \in
\Pi \times\Pi$, there exist an integer $f_{I,J}$ and scalars $k_{L,M}
\in {\mathbb C}$ for each $(L,M) \in (S \cap E_{I,J}) \cup (T \cap
E_{I,J})$ such that
\begin{equation}\label{pre-rel}
[I][J]-q^{f_{I,J}} [J][I]
= \sum_{(L,M) \in S \cap E_{I,J}} k_{L,M} [L][M]
+ \sum_{(L,M) \in T \cap E_{I,J}} k_{L,M} [L][M].
\end{equation}
To each pair $(I,J) \in \Pi\times\Pi$ and each equation as in 
(\ref{pre-rel}), we attach the integer 
$\max\{\rk L,\,(L,M) \in T \cap E_{I,J},\, k_{L,M}\neq 0\}$ if 
$\{(L,M) \in T \cap E_{I,J},\, k_{L,M}\neq 0\}
\neq\emptyset$ and $0$ otherwise.
(The definition of the rank of an element in a poset was given just after
Remark \ref{rem-before-rank}.) 
In addition, for $t \in {\mathbb N}$, we say that a pair $(I,J)$ is of
{\em type
$t$} if there exists an equation as in (\ref{pre-rel}) for which 
the attached integer is less than or equal to $t$. 
With this vocabulary, our
aim is to show that each pair $(I,J) \in \Pi\times\Pi$ is of type $0$.
We will prove the following statement: for each $t \in {\mathbb N}$, 
if $(I,J) \in \Pi\times\Pi$ is of type $t$, then $(I,J)$ is 
of type $0$; for this, we proceed by induction on $t$.
The result is trivial for $t=0$. 
Now, assume that the statement is true for pairs of type $s$ 
for some integer $s \ge 0$. 
Let $(I,J) \in \Pi\times\Pi$ be a pair of type $s+1$. 
Then there exists $f_{I,J} \in {\mathbb Z}$ and $k_{L,M} \in {\mathbb C}$ for 
$(L,M) \in (S \cap E_{I,J}) \cup (T \cap E_{I,J})$ such that
$[I][J]-q^{f_{I,J}} [J][I]
= \sum_{(L,M) \in S \cap E_{I,J}} k_{L,M} [L][M]
+ \sum_{(L,M) \in T \cap E_{I,J}} k_{L,M} [L][M]$, with in addition,
$\rk L \le s+1$ for all $(L,M) \in T \cap E_{I,J}$ such that $k_{L,M}
\neq 0$. By applying 
Proposition \ref{ASL-semipre-cr-over-complexes} 
to any pair $(L,M) \in T \cap E_{I,J}$ such that $k_{L,M} \neq 0$, 
we see that $(L,M)$ has type $s$; and so 
by the inductive hypothesis, $(L,M)$ has type $0$. 
It follows that we have an equation
$[L][M]=q^{f_{L,M}}[M][L]
+ \sum_{(L_1,M_1) \in S \cap E_{L,M}} k_{L_1,M_1} [L_1][M_1]$.
Substituting this new expression in 
$\sum_{(L,M) \in T \cap E_{I,J}} k_{L,M} [L][M]$ gives us a new expression
of $[I][J]-q^{f_{I,J}} [J][I]$ as in (\ref{pre-rel}) but 
where the second summand
in the right hand side term 
is zero, that is: $(I,J)$ is of type zero. 
This finishes the proof by induction.\qed   

\begin{subtheorem} \label{cr-in-general} -- 
Let $A$ be a commutative domain and let $u$ be any invertible element of $A$.
Let $I,J \in \Pi$; then there exists a relation
\[
[I][J]-u^{f_{I,J}} [J][I] = \sum k_{L,M} [L][M],
\]
in ${\mathcal O}_u(G_{m,n}(A))$,
where the sum is taken over pairs $(L,M)$ of elements of $\Pi$
such that $L <_\st I,J$ and $L \le_\st M$ and where $f_{I,J} \in {\mathbb Z}$
and $k_{L,M} \in A$.
\end{subtheorem}

\proof Let ${\mathcal A}={\mathbb Z}[t^{\pm 1}]$. Choose $q \in {\mathbb C}$
transcendental over ${\mathbb Q}$.  There is an injective ring homomorphism
${\mathcal A} \longrightarrow {\mathbb C}$ such that $1 \mapsto 1$ and 
$t \mapsto q$. From Remark \ref{homo-grass}, it follows that there is an
injective ring homomorphism $f \, : \, {\cal O}_t(G_{m,n}({\mathcal A}))
\longrightarrow {\cal O}_q(G_{m,n}({\mathbb C}))$ such that $1 \mapsto 1$,
$t \mapsto q$ and, for all $I \in \Pi$, $[I] \mapsto [I]$. Now, consider 
$I,J \in \Pi$~; according to Theorem \ref{ASL-cr-over-complexes}, 
there exists $f_{I,J} \in {\mathbb Z}$ and, for each $(L,M) \in \Pi$
such that $L \le_\st M$, $L <_\st I,J$, a complex number $k_{L,M}$ such that
\begin{equation} \label{r-r-r}
[I][J]-q^{f_{I,J}} [J][I] = \sum k_{L,M} [L][M],
\end{equation}
in ${\mathcal O}_q(G_{m,n}({\mathbb C}))$.
On the other hand, by Proposition \ref{basis-over-ring}, we can express the 
element $[I][J]-t^{f_{I,J}} [J][I]$ of ${\mathcal O}_t(G_{m,n}({\mathcal A}))$ 
as a linear combination of standard monomials in
${\mathcal O}_t(G_{m,n}({\mathcal A}))$:
$[I][J]-t^{f_{I,J}} [J][I]=\sum z_{L,M} [L][M]$
where the sum is taken over all pairs $(L,M)$ of elements of $\Pi$ such that
$L \le_\st M$ and $z_{L,M} \in {\mathcal A}$ for each such pair.
Now, applying $f$ to this expression gives 
$[I][J]-q^{f_{I,J}} [J][I]=\sum f(z_{L,M}) [L][M]$ 
in ${\cal O}_q(G_{m,n}({\mathbb C}))$ from which it follows,
by using relation (\ref{r-r-r}), 
the linear independence of standard monomials in 
${\mathcal O}_q(G_{m,n}({\mathbb C}))$ and the injectivity of $f$ that
$z_{L,M}$ is zero whenever $L \not <_\st I,J$.
Hence, we end up with a relation
\begin{equation} \label{comm-over-z}
[I][J]-t^{f_{I,J}} [J][I]=\sum z_{L,M} [L][M],
\end{equation}
where the sum is taken over all pairs $(L,M)$ of elements of $\Pi$ such that
$L \le_\st M$, $L <_\st I,J$ and $z_{L,M} \in {\mathcal A}$ for each such pair.
Now, there is a natural ring homomorphism ${\mathcal A} \longrightarrow A$ such
that $1 \mapsto 1$ and $t \mapsto u$. Using Remark \ref{homo-grass} again,
we get a ring homomorphism $g\, : \, {\mathcal O}_t(G_{m,n}({\mathcal A}))
\longrightarrow {\mathcal O}_u(G_{m,n}(A))$ such that 
$1 \mapsto 1$, $t \mapsto u$ and $[I] \mapsto [I]$ for all $I \in \Pi$.
It is now enough to apply $g$ to the relations (\ref{comm-over-z}) to finish
the proof.\qed\\

We are now in position to prove the main theorem of this section. 

\begin{subtheorem} -- \label{q-grass-are-qg-ASL}
Let $\Bbbk$ be a field, $q$ be a non-zero element of 
$\Bbbk$ and $m,n$ be positive integers such that $m \le n$.
Then ${\mathcal O}_q(G_{m,n}(\Bbbk))$ is a quantum graded A.S.L. on the poset
$(\Pi_{m,n},<_\st)$.
\end{subtheorem}

\proof Recall that ${\mathcal O}_q(G_{m,n}(\Bbbk))$ is an ${\mathbb N}$-graded
ring generated by the elements $[I]$, $I \in \Pi_{m,n}$ which have degree one.
By Proposition \ref{basis-over-ring}, the set of standard monomials is
a free family. It remains to use Theorems \ref{straightening-relations} and
\ref{cr-in-general} to conclude.\qed

\begin{subcorollary} -- \label{q-schub-are-qg-ASL}
Let $\Bbbk$ be a field, $q$ be a non-zero element of
$\Bbbk$ and $m,n$ be positive integers such that $m \le n$.
For $\gamma \in \Pi_{m,n}$, the quantum Schubert variety
${\mathcal O}_q(G_{m,n}(\Bbbk))_\gamma$ is a quantum 
graded A.S.L. on the poset $(\Pi_{m,n}\setminus\Omega_\gamma,<_\st)$, where 
$\Omega_\gamma=\{\pi \in \Pi_{m,n} \tq \pi \not\ge_\st \gamma\}$.
\end{subcorollary}

\proof For $\gamma \in \Pi_{m,n}$, ${\mathcal O}_q(G_{m,n}(\Bbbk))_\gamma=
{\mathcal O}_q(G_{m,n}(\Bbbk))/\langle\Omega_\gamma\rangle$. Since
$\Omega_\gamma$ is a $\Pi_{m,n}$-ideal, the result follows from 
Theorem \ref{q-grass-are-qg-ASL} and Corollary \ref{quotients-graded-ASL}.\qed

\subsection{Quantum determinantal rings are quantum graded A.S.L.}

Throughout this subsection, $\Bbbk$ denotes a field.\\

In this subsection, we want to derive the fact that quantum determinantal rings 
are quantum graded A.S.L. from Theorem \ref{q-grass-are-qg-ASL}. 
For this, we need to relate the ring of quantum matrices and quantum 
grassmannians, and this is done using a {\it dehomogenisation map} from 
[KLR] the construction of which we now recall.\\

From [KLR; Lemma 1.5], we know that $[n+1,\dots,n+m]$ is a normal element 
in ${\mathcal O}_q(G_{m,m+n}(\Bbbk))$. Hence, we may localise 
${\mathcal O}_q(G_{m,m+n}(\Bbbk))$ with respect to the set of powers of 
$[n+1,\dots,n+m]$; this localisation is denoted 
${\mathcal O}_q(G_{m,m+n}(\Bbbk))[[n+1,\dots,n+m]^{-1}]$. For convenience, if 
$I=\{i_1 < \dots < i_m\}$ with $1 \le i_1 < \dots < i_m \le m+n$, we put
$\{I\}=[I][n+1,\dots,n+m]^{-1} 
\in {\mathcal O}_q(G_{m,m+n}(\Bbbk))[[n+1,\dots,n+m]^{-1}]$.
Now, let $\phi$ be the automorphism of ${\mathcal O}_q(M_{m,n}(\Bbbk))$
defined by $\phi(X_{ij})=q^{-1}X_{ij}$ for $1 \le i \le m$ and $1 \le j \le n$.
From [KLR; Corollary 4.4], we know that there is an isomorphism of 
$\Bbbk$-algebras
\[
\begin{array}{ccrcl}
D_{m,n} & : & {\mathcal O}_q(M_{m,n}(\Bbbk))[y,y^{-1};\phi] & \longrightarrow
& {\mathcal O}_q(G_{m,m+n}(\Bbbk))[[n+1,\dots,n+m]^{-1}] \cr 
\end{array}
\] 
defined by 
$D_{m,n}(X_{ij})=\{\{j,n+1,\dots,\widehat{n+m+1-i},\dots,n+m\}\}$
and $D_{m,n}(y)=[n+1,\dots,n+m]$.
(Here $\{j,n+1,\dots,\widehat{n+m+1-i},\dots,n+m\}$ stands for the (ordered)
set obtained from $\{j,n+1,\dots,n+m\}$ by removing $n+m+1-i$.)\\

Recall from Subsection \ref{subsec-q-grass} the definition of $\Delta_{m,n}$.  
In addition, to any pair $(I,J) \in \Delta_{m,n}$, 
with $I=\{i_1 < \dots < i_t\}$, $J=\{j_1 < \dots < j_t\}$, $1 \le t \le m$,
we associate the index set 
$K_{(I,J)} \in \Pi_{m,n+m}$ obtained by ordering the elements of the set
$\{j_1,\dots,j_t,n+1,\dots,n+m\}\setminus\{n+m+1-i_1,\dots,n+m+1-i_t\}$.

\begin{sublemma} -- \label{the-map-Dehom}
We keep the above notation. Let $(I,J) \in \Delta_{m,n}$ with
$I=\{i_1< \dots <i_t\} \subseteq \{1,\dots,m\}$,
$J=\{j_1< \dots <j_t\} \subseteq \{1,\dots,n\}$ and $1 \le t \le m$. Then
\[
D_{m,n}([I|J])=\{K_{(I,J)}\}=[K_{(I,J)}][\{n+1,\dots,n+m\}]^{-1}.
\]
\end{sublemma}

\proof The proof is by induction on $t$, the case $t=1$ being nothing but the
definition of $D_{m,n}$. Now, let $t$ be an integer such that $1 \le t < m$ and
assume that 
the result is true for $t \times t$ minors. We consider $(I,J) \in 
\Delta_{m,n}$ with
$I=\{i_1< \dots <i_{t+1}\} \subseteq \{1,\dots,m\}$,
$J=\{j_1< \dots <j_{t+1}\} \subseteq \{1,\dots,n\}$ and, for $1 \le k \le t+1$,
we let $I_k$ (respectively $J_k$) be the index set obtained from $I$ 
(respectively $J$)
by removing the $k$-th entry. Then, by [PW; Corollary 4.4.4], we have
\[
[I|J]= \sum_{k=1}^{t+1} (-q)^{t+1-k} [I_{t+1}|J_k]X_{i_{t+1},j_k}.
\]
Hence, by using the inductive hypothesis and [KLR; Lemma 1.5], we get
\[
\begin{array}{rcl}
D_{m,n}([I|J])
&=& \displaystyle\sum_{k=1}^{t+1} 
(-q)^{t+1-k} \{K_{(I_{t+1},J_k)}\}\{K_{(i_{t+1},j_k)}\} \cr
& & \cr
&=& q\left(
\displaystyle\sum_{k=1}^{t+1} 
(-q)^{t+1-k} [K_{(I_{t+1},J_k)}][K_{(i_{t+1},j_k)}]\right)
[n+1,\dots,n+m]^{-2}.
\end{array}
\]
Now, by applying [KLR; Theorem 2.5] in 
${\mathcal O}_q(G_{m,m+n}(\Bbbk))$ with
$K=(J \cup \{n+1,\dots,n+m\})\setminus\{m+n+1-i_1,\dots,m+n+1-i_t\}$,
$J_1=\emptyset$ and 
$J_2=\{n+1,\dots,n+m\}\setminus\{n+m+1-i_{t+1}\}$ we obtain the quantum 
Pl\"ucker relation 
\[
\sum_{k=1}^{t+1} (-q)^{m-(k-1)}[K_{(I_{t+1},J_k)}][K_{(i_{t+1},j_k)}]
+ (-q)^{m-t-1} [K_{(I,J)}][n+1,\dots,n+m]=0.
\]
(Notice that the first summand corresponds to the case where $K'' \subseteq J$ 
while the second corresponds to the case where $K''=\{m+n+1-i_{t+1}\}$.)
We may rewrite this last relation as
\[
q\sum_{k=1}^{t+1} (-q)^{t+1-k}
[K_{(I_{t+1},J_k)}][K_{(i_{t+1},j_k)}]
=[K_{(I,J)}][n+1,\dots,n+m].
\]
Hence, we get
\[
D_{m,n}([I|J])=[K_{(I,J)}][n+1,\dots,n+m]^{-1}=\{K_{(I,J)}\}.
\]
The proof is complete.\qed\\

On the set $\Delta_{m,n}$ introduced above, we put a partial order that 
we also denote by 
$\le_\st$. Let $u,v$ be integers such that $1 \le u,v \le m$ and let 
$(I,J)$ and
$(K,L)$ be an index pairs with
$I=\{i_1< \dots <i_u\} \subseteq \{1,\dots,m\}$,
$J=\{j_1< \dots <j_u\} \subseteq \{1,\dots,n\}$,
$K=\{k_1< \dots <k_v\} \subseteq \{1,\dots,m\}$,
$L=\{l_1< \dots <l_v\} \subseteq \{1,\dots,n\}$.
We define $\le_\st$ as follows:
\[
(I,J) \le_\st (K,L) 
\Longleftrightarrow 
\left\{
\begin{array}{l}
u \ge v, \cr 
i_s \le k_s \quad\mbox{for}\quad 1 \le s \le v , \cr
j_s \le l_s \quad\mbox{for}\quad 1 \le s \le v .
\end{array}
\right.
\]
The following combinatorial lemma relates
the orders on $\Delta_{m,n}$ and $\Pi_{m,n+m}$. Notice first that the map
\[
\begin{array}{ccrcl}
\delta_{m,n} & : & \Delta_{m,n} & \longrightarrow & \Pi_{m,n+m}
\setminus\{\{n+1<\dots<n+m\}\}
\cr
 & & (I,J) & \mapsto & K_{(I,J)}
\end{array}
\]
is bijective.

\begin{sublemma} -- \label{resp-order}
For $(I,J),(K,L) \in \Delta_{m,n}$, one has $(I,J) \le_\st (K,L)$
if and only if $\delta_{m,n}((I,J)) \le_\st \delta_{m,n}((K,L))$. 
Hence, $\delta_{m,n}$ is a bijective
map of posets between $(\Delta_{m,n},\le_\st)$ and 
$(\Pi_{m,n+m}\setminus\{\{n+1<\dots<n+m\}\},\le_\st)$.
\end{sublemma}

\proof See the proof of [BV; Lemma 4.9].\qed

\begin{subtheorem} \label{mat-q-q-gr-ASL} -- 
Let $\Bbbk$ be a field, $m,n$ be positive integers such that $m \le n$
and $q$ be any element of $\Bbbk^\ast$. Then
${\mathcal O}_q(M_{m,n}(\Bbbk))$ is a quantum graded A.S.L. on the poset
$(\Delta_{m,n},\le_\st)$.
\end{subtheorem}

\proof Clearly, ${\mathcal O}_q(M_{m,n}(\Bbbk))$ is an ${\mathbb N}$-graded 
$\Bbbk$-algebra 
and the elements of $\Delta_{m,n}$ are homogeneous with positive degree
and generate ${\mathcal O}_q(M_{m,n}(\Bbbk))$ as a $\Bbbk$-algebra.
Recall, from Lemma \ref{the-map-Dehom}, that 
$D_{m,n}([I|J])=[\delta_{m,n}((I,J))][n+1,\dots,n+m]^{-1}$, 
for any index pair $(I,J) \in \Delta_{m,n}$. 
As a consequence, if
$(I_1,J_1) \le_\st \dots \le_\st (I_t,J_t)$ is an 
increasing sequence of elements of $\Delta_{m,n}$, 
then there exists $\alpha \in {\mathbb Z}$
such that $D_{m,n}([I_1|J_1] \dots [I_t|J_t])
=q^\alpha[\delta_{m,n}((I_1,J_1))]\dots [\delta_{m,n}((I_t,J_t))]
[n+1,\dots,n+m]^{-t}$  
(here we use [KLR; Lemma 1.5]). From 
this observation, and using Lemma \ref{resp-order} 
and the linear independence of
standard monomials in ${\mathcal O}_q(G_{m,m+n}(\Bbbk))$, it follows 
easily that 
the set of standard monomials in ${\mathcal O}_q(M_{m,n}(\Bbbk))$ form a free 
family. Now, consider two incomparable elements $(I,J)$ and $(K,L)$ of 
$\Delta_{m,n}$. 
By Lemma \ref{resp-order}, $\delta_{m,n}((I,J))$ and $\delta_{m,n}((K,L))$ are 
incomparable elements of $\Pi_{m,n+m}$. 
Then, by Theorem \ref{straightening-relations},
there exist $K_i,L_i \in \Pi_{m,n+m}$ and 
$k_i \in \Bbbk$, for
$1 \le i \le t$, such that 
\[
[\delta_{m,n}((I,J))][\delta_{m,n}((K,L))]=\sum_{i=1}^t k_i [K_i][L_i] 
\]
with, in addition $K_i \le_\st L_i$ and $K_i <_\st \delta_{m,n}((I,J)),
\delta_{m,n}((K,L))$, for
$1 \le i \le t$. If we put $M=\{n+1,\dots,n+m\} \in \Pi_{m,n+m}$, 
by [KLR; Lemma 1.5],
we then have 
\[
[\delta_{m,n}((I,J))][M]^{-1}[\delta_{m,n}((K,L))][M]^{-1}
=\sum_{i=1}^t k_i^\prime [K_i][M]^{-1}[L_i][M]^{-1},
\]
where $k_i^\prime \in \Bbbk$, for $1 \le i \le t$. At this point, 
it is worth 
mentioning that no $K_i$ can equal $M$ since 
$K_i <_\st \delta_{m,n}((I,J)),\delta_{m,n}((K,L))$ and 
$M$ is maximal in $\Pi_{m,n+m}$. 
In contrast, it is
possible for an $L_i$ to equal $M$; so that in the above equation, a term
$[L_i][M]^{-1}$ might very well equal $1$. Now, by 
applying $D_{m,n}^{-1}$ to this
last equation we obtain an expression for  $[I|J][K|L]$ as a linear
combination of terms of the form $[E|F][G|H]$ or $[E|F]$, where $(E,F),(G,H) 
\in \Delta_{m,n}$, 
$(E,F) \le_\st (G,H)$ and $(E,F) <_\st (I,J),(K,L)$. This shows that
condition (4) in the definition of a quantum graded A.S.L. is satisfied.
Finally, consider two elements $(I,J)$ and $(K,L)$ in $\Delta_{m,n}$. By 
Theorem \ref{cr-in-general} 
there exist $f \in {\mathbb Z}$, elements  $K_i,L_i \in 
\Pi_{m,n+m}$ and 
$k_i \in \Bbbk$, for $1 \le i \le t$, such that
\[
[\delta_{m,n}((I,J))][\delta_{m,n}((K,L))]
-q^f[\delta_{m,n}((K,L))][\delta_{m,n}((I,J))]
=\sum_{i=1}^t k_i [K_i][L_i]
\]
with, in addition, $K_i \le_\st L_i$ and $K_i <_\st 
\delta_{m,n}((I,J)),\delta_{m,n}((K,L))$, for
$1 \le i \le t$. Now, applying the same procedure as above, we end up with
an expression for $[I|J][K|L]-q^g[K|L][I|J]$ (for some $g \in {\mathbb Z}$)
as a linear combination of terms of the form $[E|F][G|H]$ or $[E|F]$, 
where $(E,F),(G,H) \in \Delta_{m,n}$ with  
$(E,F) \le_\st (G,H)$ and $(E,F) <_\st (I,J),(K,L)$. This shows that
condition (5) in the definition of a quantum graded A.S.L. is satisfied.\qed\\

For $1 \le t \le m$, we denote by ${\mathcal I}_t$ the ideal of 
${\mathcal O}_q(M_{m,n}(\Bbbk))$ generated by the $t \times t$ quantum minors. 
On the other hand, for $t \ge 2$, 
we denote by $\Omega_t$ the $\Delta_{m,n}$-ideal defined by
\[
\Omega_t=\{(I,J) \in \Delta_{m,n} 
\tq (I,J) \not\ge_\st (\{1,\dots,t-1\},\{1,\dots,t-1\})\},
\]
and we put $\Omega_1=\Delta_{m,n}$. 
For $1 \le t \le m$, 
it is easy to check that $\Omega_t$ is just the set
of index pairs $(I,J)$ such that $t \le |I|=|J| \le m$. 
Hence, 
${\mathcal I}_t=\langle\Omega_t\rangle$, by 
using the well known Laplace expansions
for quantum minors (see [PW; Corollary 4.4.4]).

\begin{subcorollary} \label{q-det-q-gr-asl} -- 
Let $\Bbbk$ be a field, $m,n$ be positive integers such that $m \le n$ and
$q$ be any element of $\Bbbk^\ast$.
For $1 \le t \le m$, the quantum determinantal ring 
${\mathcal O}_q(M_{m,n}(\Bbbk))/{\mathcal I}_t$ is a quantum graded A.S.L. 
on the poset $\Delta_{m,n} \setminus \Omega_t$.
\end{subcorollary}

\proof According to the comments above, this is an immediate consequence of
Theorem \ref{mat-q-q-gr-ASL} and Corollary \ref{quotients-graded-ASL}.\qed

\section{Quantum determinantal rings and quantum grassmannians are 
AS-Cohen-Macaulay.}

In this section, we reach our main aim. We prove that quantum determinantal
rings and quantum Grassmannians are AS-Cohen-Macaulay. Once this is done,
we determine which of these rings are 
AS-Gorenstein, by means of their Hilbert series. 

\begin{theorem} -- \label{q-detrings-are-ASCM}
Let $\Bbbk$ be a field, $m,n$ be positive integers such that $m \le n$ and
$q$ be any element of $\Bbbk^\ast$.
For $1 \le t \le m$, the quantum determinantal ring
${\mathcal O}_q(M_{m,n}(\Bbbk))/{\mathcal I}_t$ is AS-Cohen-Macaulay.
\end{theorem}

\proof The case $t=1$ is trivial, thus we assume that $t \ge 2$.
First notice that $\Delta_{m,n}\setminus\Omega_t=
\{(I,J) \in \Delta_{m,n} \tq (I,J) 
\ge_\st (\{1,\dots,t-1\},\{1,\dots,t-1\})\}$ which is
a distributive lattice by [BV; Theorem 5.3] and hence a wonderful poset 
(see [BV; p. 58]). Thus, the result follows from 
Theorem \ref{critirion-ASL-ASCM} in conjunction with
Corollary \ref{q-det-q-gr-asl}. \qed

\begin{theorem} -- \label{q-grass-schub-are-ASCM}
Let $\Bbbk$ be a field, $q$ be a non-zero element of
$\Bbbk$ and $m,n$ be positive integers such that $m \le n$.
For $\gamma \in \Pi_{m,n}$, the quantum Schubert variety
${\mathcal O}_q(G_{m,n}(\Bbbk))_\gamma$ is AS-Cohen-Macaulay.
In particular, the quantum Grassmannian
${\mathcal O}_q(G_{m,n}(\Bbbk))$ is AS-Cohen-Macaulay.
\end{theorem}

\proof Recall that $\Pi_{m,n}$ has a single minimal element. Hence the
quantum grassmannian is a special case of a quantum Schubert variety.
Let $\gamma \in \Pi_{m,n}$. From Corollary \ref{q-schub-are-qg-ASL}, 
the quantum Schubert variety
${\mathcal O}_q(G_{m,n}(\Bbbk))_\gamma$ is a quantum
graded A.S.L. on the poset $(\Pi_{m,n}\setminus\Omega_\gamma,<_\st)$, where
$\Omega_\gamma=\{\pi \in \Pi_{m,n} \tq \pi \not\ge_\st \gamma\}$. On the other
hand, by [BV; Theorem 5.4], 
the poset $(\Pi_{m,n}\setminus\Omega_\gamma,<_\st)$ is a distributive
lattice and hence a wonderful poset
(see [BV; p. 58]). Thus, the result follows from
Theorem \ref{critirion-ASL-ASCM}. \qed\\

We now show that quantum grassmannians are AS-Gorenstein and determine which
quantum determinantal rings are AS-Gorenstein. This can be easily deduced
from the above results using the criterion given in [JZ; Theorem 6.2].\\

Recall that, if $A=\oplus_{i \in {\mathbb N}} A_i$ is a noetherian 
${\mathbb N}$-graded connected $\Bbbk$-algebra, then it is locally finite and we
can speak of its Hilbert series, namely
\[
H_A(t)=\sum_{i \ge 0} (\dim_\Bbbk A_i)t^i \in {\mathbb Q}[[t]].
\]

If, in addition, $A$ is AS-Cohen-Macaulay, a domain and has enough normal 
elements, then [JZ; Theorem 6.2] shows that the fact that $A$ is AS-Gorenstein
or not can be read off from its Hilbert series. This is the criterion we will
use to prove the following theorems.

\begin{remark} \label{lien-ASgor-gor} -- \rm
Let $R$ be a noetherian commutative ring. Then, $R$ is
said to be Gorenstein if, for every maximal ideal ${\mathfrak p}$ of $R$,
the local ring $R_{\mathfrak p}$ is of finite self-injective dimension. 
If, moreover, the Krull dimension of $R$ is finite, then $R$ is Gorenstein if
and only if it has finite self-injective dimension (see 
[B; Theorem and definition \S1]).
Now, let $A=\oplus_{i \in {\mathbb N}} A_i$ be a commutative
noetherian ${\mathbb N}$-graded connected $\Bbbk$-algebra. Since, clearly,
$A$ has finite Krull dimension, $A$ is Gorenstein if and only if
$A$ has finite (self-)injective dimension. On the other hand it is clear that 
$A$ has enough normal elements, since $A$ is commutative. Thus, from Remark 
\ref{ref-co-case}, it follows that $A$ is AS-Gorenstein if and only if it is
Gorenstein.
\end{remark} 

\begin{theorem} -- \label{q-grass-are-ASG} Let $\Bbbk$ be a field, $q$
be a non-zero element of $\Bbbk$ and $m,n$ be positive integers
such that $m \le n$.  Then
the quantum Grassmannian ${\mathcal O}_q(G_{m,n}(\Bbbk))$ is
AS-Gorenstein.  \end{theorem}

\proof The quantum grassmannian ${\mathcal O}_q(G_{m,n}(\Bbbk))$ is a
domain since it is a subalgebra of the domain ${\mathcal
O}_q(M_{m,n}(\Bbbk))$.  It is AS-Cohen-Macaulay by Theorem
\ref{q-grass-schub-are-ASCM}.  And, it has enough normal elements since
it is a quantum graded A.S.L.  (see Theorem \ref{q-grass-are-qg-ASL} and
Remark \ref{ASL-ENE}). 

On the other hand, it follows at once from Proposition 
\ref{basis-over-ring} that for all $q \in \Bbbk^\ast$, the Hilbert
series of ${\mathcal O}(G_{m,n}(\Bbbk)):={\mathcal O}_1(G_{m,n}(\Bbbk))$
and ${\mathcal O}_q(G_{m,n}(\Bbbk))$ coincide.  Hence, using [JZ;
Theorem 6.2], it is enough to prove that the (usual) homogeneous
coordinate ring ${\mathcal O}(G_{m,n}(\Bbbk))$ is AS-Gorenstein. 
However, this follows at once from Remark \ref{lien-ASgor-gor} and [BV;
Corollary 8.13].\qed

\begin{theorem} -- \label{q-detrings-and-ASG}
Let $\Bbbk$ be a field, $m,n$ be positive integers such that $m \le n$ and
$q$ be any element of $\Bbbk^\ast$.
For $1 \le t \le m$, the quantum determinantal ring
${\mathcal O}_q(M_{m,n}(\Bbbk))/{\mathcal I}_t$ is AS-Gorenstein if and only
if $t=1$ or $m=n$.
\end{theorem}

\proof The quantum determinantal rings ${\mathcal
O}_q(M_{m,n}(\Bbbk))/{\mathcal I}_t$ are domains by [GL; Corollary 2.6]. 
They are AS-Cohen-Macaulay by Theorem \ref{q-detrings-are-ASCM}, and 
they have enough normal elements since they are quantum graded A.S.L. 
(see Corollary \ref{q-det-q-gr-asl} and Remark \ref{ASL-ENE}). 

On the other hand, it follows at once from Theorem \ref{q-det-q-gr-asl}
and Proposition \ref{standard-mon-basis}
that, for all $q \in \Bbbk^\ast$, the Hilbert series of ${\mathcal
O}(M_{m,n}(\Bbbk))/{\mathcal I}_t
:={\mathcal O}_1(M_{m,n}(\Bbbk))/{\mathcal I}_t$ and ${\mathcal
O}_q(M_{m,n}(\Bbbk))/{\mathcal I}_t$ coincide.  Hence, using [JZ; Theorem 6.2],
it is
enough to prove that the (usual) coordinate ring
${\mathcal O}(M_{m,n}(\Bbbk))/{\mathcal I}_t$ is AS-Gorenstein
if and only if $t=1$ or $m=n$. However, this follows at once
from Remark \ref{lien-ASgor-gor} and [BV; Corollary 8.9].\qed


\section*{References.}
{\bf [ATV]} M. Artin, J. Tate and M. Van den Bergh. {\em Some algebras
associated to automorphisms of elliptic curves}, The Grothendieck
Festschrift, Vol.  I, 33-85, Progr.  Math., {\bf 86}, Birkh\"auser
Boston, Boston, MA, 1990.
\newline
{\bf [B]} H. Bass. {\em On the ubiquity of Gorenstein rings}. Math. Zeitschr.
{\bf 82} (1963), 8--28. \newline
{\bf [BH]}  W. Bruns and J. Herzog. Cohen-Macaulay rings. 
Cambridge Studies in Advanced Mathematics, 39.
Cambridge University Press, Cambridge, 1998.\newline
{\bf [BV]} W. Bruns and U. Vetter. Determinantal rings.  
Lecture notes in Mathematics, 1327. Springer-Verlag, Berlin,
1988. \newline
{\bf [DEP]} C. De Concini, D. Eisenbud and C. Procesi. Hodge algebras. 
Ast\'erisque, {\bf 91}. Soci\'et\'e Math\'ematique de France, Paris, 1982.
\newline
{\bf [E]} D. Eisenbud. {\em Introduction to algebras with straightening laws.}
Ring theory and algebra, III, 243--268. Lecture Notes in Pure and Appl. Math., 
55, Dekker, New-York, 1980.
\newline
{\bf [GL]} K.R. Goodearl and T.H. Lenagan. {\em Quantum determinantal ideals.}
Duke Math. J. {\bf 103} (2000), 165-190.   \newline
{\bf [GLR]}  K. R. Goodearl, T. H. Lenagan and L. Rigal. {\em The first
fundamental theorem of coinvariant theory for the quantum general linear
group}. Publ. RIMS (Kyoto) {\bf 36} (2000), 269-296.
\newline
{\bf [HL(1)]} T. J. Hodges and T. Levasseur. 
{\em Primitive ideals of ${\mathbb C}_q[SL(3)]$}. Comm. Math. Physics
{\bf 156} (1993), 581-605.
\newline
{\bf [HL(2)]} T. J. Hodges and T. Levasseur.
{\em  Primitive ideals of ${\mathbb C}_q[SL(n)]$}. J. of Algebra {\bf 168}
(1994), 455-468.
\newline
{\bf [JZ]} P. J\o rgensen and J.J. Zhang. 
{\em Gourmet's guide to Gorensteinness}.  Adv. Math.  {\bf 151}  (2000),  
313--345.
\newline  
{\bf [KLR]} A. C. Kelly, T.H. Lenagan and L. Rigal. {\em  Ring
theoretic properties of quantum grassmannians}. 
To appear in Journal of Algebra and its Applications 
\newline 
{\bf [KL]} G. Krause and T.H. Lenagan.  Growth of algebras and
Gelfand-Kirillov dimension, Revised edition.
Graduate Studies in Mathematics, 22. American Mathematical Society,
Providence, RI, 2000.\newline
{\bf [LR]} T.H. Lenagan and L. Rigal.   
{\em The maximal order property for
quantum determinantal rings}. 
Proc. Edin. Math. Soc. {\bf 46} (2003), 513-529. \newline
{\bf [Lev]} T. Levasseur. 
{\em Some properties of non-commutative regular graded
rings}. Glasgow Math. J. {\bf 34} (1992), 277-300. \newline
{\bf [PW]} B. Parshall and J.-P. Wang.  {\em Quantum linear groups}. Mem.
Amer. Math. Soc {\bf 89} (1991), no. 439.
\newline
{\bf [W]} C.A. Weibel. 
An introduction to homological algebra. Cambridge Studies in
Advanced Mathematics, 38. Cambridge University Press, Cambridge, 1994.
\newline
{\bf [YZ]} A. Yekutieli and J.J. Zhang. {\em Rings with Auslander dualizing
complexes}. Journal of Algebra {\bf 213} (1999), 1-51.
\newline
{\bf [Z]} J.J. Zhang. {\em Connected graded Gorenstein algebras with enough
normal elements}. Journal of Algebra {\bf 189} (1997), 390-405. 
\newline 


\vskip 1cm

\noindent T. H. Lenagan: \\
School of Mathematics, University of Edinburgh,\\
James Clerk Maxwell Building, King's Buildings, Mayfield Road,\\
Edinburgh EH9 3JZ, Scotland\\
E-mail: tom@maths.ed.ac.uk \\
\\
L. Rigal: \\
Universit\'e Jean Monnet
(Saint-\'Etienne), \\
Facult\'e des Sciences et
Techniques, \\
D\'e\-par\-te\-ment de Math\'ematiques,\\
23 rue du Docteur Paul Michelon,\\
42023 Saint-\'Etienne C\'edex 2,\\France\\
E-mail: Laurent.Rigal@univ-st-etienne.fr


\end{document}